\documentclass[12pt]{article}
\usepackage{amsfonts,amssymb,amsmath}
\newtheorem{theorem}{Theorem}[section]
\newtheorem{proposition}{Proposition}[section]
\newtheorem{lemma}{Lemma}[section]
\newtheorem{corollary}{Corollary}[section]
\newtheorem{definition}{Definition}[section]
\newtheorem{remark}{Remark}[section]

\begin{document}

\title{Extension of Finite Rank Operators and Local Structures in Operator Ideals}
\author{Frank Oertel \\
University of Bonn\\
Department of Statistics\\
Adenauerallee 24--26\\
D - 53113 Bonn\\
(e-mail: {\em oertel@addi.finasto.uni--bonn.de})}
\maketitle

\begin{abstract}
{\noindent}We develop general techniques and present an approach to solve
the problem of constructing a maximal Banach ideal $({\frak A},{\bf A)}$
which does {\it not} satisfy a transfer of the norm estimation in the
principle of local reflexivity to its norm ${\bf A}$. This approach leads us
to the investigation of product operator ideals containing ${\frak L}_2$
(the collection of all Hilbertian operators) as a factor. Using the local
properties of such operator ideals -- which are typical examples of ideals
with property (I) and property (S) --, trace duality and an extension of
suitable finite rank operators even enable us to show that ${\frak L}_\infty 
$ cannot be totally accessible -- answering an open question of Defant and
Floret.\newline

{\noindent}{\it {Key words and phrases:}} Accessibility, Banach spaces,
cotype 2, finite rank operators, Hilbert space factorization, Grothendieck's
Theorem, operator ideals, principle of local reflexivity, tensor norms%
\newline

{\noindent}{\it {1991 AMS Mathematics Subject Classification:}} primary
46M05, 47D50; secondary 47A80.
\end{abstract}

\section{Introduction}

The aim of the present paper is to present a thorough investigation of
operator ideals $({\frak A},{\bf A})$ in relation to a transfer of the norm
estimation in the classical principle of local reflexivity to their ideal
(quasi--)norm ${\bf A}$. In particular, we are interested in constructing
examples of maximal Banach ideals which do not satisfy such a transfer. Due
to the local nature of this {\it principle of local reflexivity for operator
ideals} (called ${\frak A}-LRP$) -- which had been introduced and discussed
in \cite{oe1} and \cite{oe2} -- and the local nature of maximal Banach
ideals, local versions of injectivity (right--accessibility) resp.
surjectivity (left--accessibility) of suitable operator ideals and
factorizations through operators with finite dimensional range even imply
interesting relations between operators with infinite dimensional range.
After {\it extending} finite rank operators in certain quasi--Banach ideals $%
{\frak A}$, the ${\frak A}-LRP$ and the calculation of conjugate ideal norms
then allow us to neglect the structure of the range space, so that we may
leave the finite dimensional case. We will give sufficient conditions on $%
{\frak A}$ to guarantee that each finite rank operator $L$ has a {\it finite
rank--extension} $\widetilde{L}$ so that ${\bf A}(\widetilde{L})\leq
(1+\epsilon )\cdot {\bf A}(L)$ -- for given $\epsilon >0$. Consequently, we
are lead to the problem under which circumstances a finite rank operator $%
L\in {\frak A}\circ {\frak B}$ has a factorization $L=AB$ so that ${\bf A}%
(A)\cdot {\bf B}(B)\leq (1+\epsilon )\cdot {\bf A}\circ {\bf B}(L)$ and $A$
resp. $B$ has finite dimensional range. Operator ideals ${\frak A}\circ 
{\frak B}$ with such a {\it property (I)} resp. {\it property (S)} had been
introduced in \cite{jo} to prepare a detailed investigation of trace ideals.

After introducing the necessary framework which also includes a full
description of the technical concept of ultrastability, we recall the
definition of the ${\frak A}-LRP$ and its first consequences. Not only in
view of looking for a counterexample of a maximal Banach ideal $({\frak A}_0,%
{\bf A}_0)$ which does not satisfy the ${\frak A}_0-LRP$, we will see that
the property (I) of ${\frak A}^{*}\circ {\frak L}_\infty $ plays a
fundamental part in this paper; it even enables us to show that ${\frak L}%
_\infty $ is not totally accessible -- answering a question of Defant and
Floret (see Theorem 4.1)! We finish the paper with further applications,
where we also consider linear operators acting between Banach spaces with
cotype 2 which do not have the approximation property (such as Pisier`s
space $P$). We apply the machinery of section 3 to product operator ideals
which contain the operator ideal $({\frak L}_2,{\bf L}_2)$ as a factor and
reveal surprising relations between the principle of local reflexivity for
the maximal hull of such operator ideals and the existence of an ideal--{\it %
norm} on these product ideals.

\section{The framework}

In this section, we introduce the basic notation and terminology which we
will use throughout in this paper. We only deal with Banach spaces and most
of our notations and definitions concerning Banach spaces and operator
ideals are standard. We refer the reader to the monographs \cite{df}, \cite
{djt} and \cite{p1} for the necessary background in operator ideal theory
and the related terminology. Infinite dimensional Banach spaces over the
field ${\Bbb K\in \{R},{\Bbb C\}}$ are denoted throughout by $W,X,Y$ and $Z$
in contrast to the letters $E,F$ and $G$ which are used for finite
dimensional Banach spaces only. The space of all operators (continuous
linear maps) from $X$ to $Y$ is denoted by ${\frak L}(X,Y)$, and for the
identity operator on $X$, we write $Id_X$. The collection of all finite rank
(resp. approximable) operators from $X$ to $Y$ is denoted by ${\frak F}(X,Y)$
(resp. $\overline{{\frak F}}(X,Y)$), and ${\frak E}(X,Y)$ indicates the
collection of all operators, acting between finite dimensional Banach spaces 
$X$ and $Y$ (elementary operators). The dual of a Banach space $X$ is
denoted by $X^{\prime }$, and $X^{\prime \prime }$ denotes its bidual $%
(X^{\prime })^{\prime }$. If $T\in {\frak L}(X,Y)$ is an operator, we
indicate that it is a metric injection by writing $T:X\stackrel{1}{%
\hookrightarrow }Y$, and if it is a metric surjection, we write $T:X%
\stackrel{1}{\twoheadrightarrow }Y$. If $X$ is a Banach space, $E$ a finite
dimensional subspace of $X$ and $K$ a finite codimensional subspace of $X$,
then $B_X:=\{x\in X:\Vert x\Vert \leq 1\}$ denotes the closed unit ball, $%
J_E^X$ $:E\stackrel{1}{\hookrightarrow }X$ the canonical metric injection
and $Q_K^X:X\stackrel{1}{\twoheadrightarrow }X\diagup K$ the canonical
metric surjection. Finally, $T^{\prime }\in {\frak L}(Y^{\prime },X^{\prime
})$ denotes the dual operator of $T\in {\frak L}(X,Y)$.

If $({\frak A},{\bf A})$ and $({\frak B},{\bf B})$ are given quasi--Banach
ideals, we will use throughout the shorter notation $({\frak A}^d,{\bf A}^d)$
for the dual ideal and the abbreviation ${\frak A}\stackrel{1}{=}{\frak B}$
for the isometric equality $({\frak A},{\bf A})=({\frak B},{\bf B})$. We
write ${\frak A}\subseteq {\frak B}$ if, regardless of the Banach spaces $X$
and $Y$, we have ${\frak A}(X,Y)\subseteq {\frak B}(X,Y)$. If $X_0$ is a
fixed Banach space, we write ${\frak A}(X_0,\cdot )\subseteq {\frak B}%
(X_0,\cdot )$ (resp. ${\frak A}(\cdot ,X_0)\subseteq {\frak B}(\cdot ,X_0)$)
if, regardless of the Banach space $Z$ we have ${\frak A}(X_0,Z)\subseteq 
{\frak B}(X_0,Z)$ (resp. ${\frak A}(Z,X_0)\subseteq {\frak B}(Z,X_0)$). The
metric inclusion $({\frak A},{\bf A})\subseteq ({\frak B},{\bf B})$ is often
shortened by ${\frak A}\stackrel{1}{\subseteq }{\frak B}$. If ${\bf B(}%
T)\leq {\bf A}(T)$ for all finite rank (resp. elementary) operators $T\in 
{\frak F}$ (resp. $T\in {\frak E}$), we sometimes use the abbreviation $%
{\frak A}\stackrel{\frak F}{\subseteq }{\frak B}$ (resp. ${\frak A}\stackrel%
{\frak E}{\subseteq }{\frak B}$).

First we recall the basic notions of Grothendieck's metric theory of tensor
products (cf., eg., \cite{df}, \cite{gl}, \cite{gr}, \cite{l}), which
together with Pietsch's theory of operator ideals spans the mathematical
frame of this paper. A {\it {tensor norm}} $\alpha $ is a mapping which
assigns to each pair $(X,Y)$ of Banach spaces a norm $\alpha (\cdot ;X,Y)$
on the algebraic tensor product $X\otimes Y$ (shorthand: $X$ ${\otimes }%
_\alpha $ $Y$ and $X\tilde{\otimes}_\alpha Y$ for the completion) so that

\begin{itemize}
\item  $\varepsilon \leq \alpha \leq \pi $

\item  $\alpha $ satisfies the metric mapping property: If $S\in {\cal L}%
(X,Z)$ and $T\in {\cal L}(Y,W)$, then $\Vert S\otimes T:X\otimes _\alpha
Y\longrightarrow Z\otimes _\alpha W\Vert \leq \Vert S\Vert \ \Vert T\Vert $ .
\end{itemize}

\noindent Wellknown examples are the injective tensor norm $\varepsilon $,
which is the smallest one, and the projective tensor norm $\pi $, which is
the largest one. For other important examples we refer to \cite{df}, \cite
{gl}, or \cite{l}. Each tensor norm $\alpha $ can be extended in two natural
ways. For this, denote for given Banach spaces $X$ and $Y$ 
\[
\text{FIN}(X):=\{E\subseteq X\mid E\in \text{FIN}\}\hspace{0.2cm}\text{and}%
\hspace{0.2cm}\text{COFIN}(X):=\{L\subseteq X\mid X/L\in \text{FIN}\}\text{,}
\]
where FIN stands for the class of all finite dimensional Banach spaces. Let $%
z\in X\otimes Y$. Then the {\it {finite hull}}\/ $\stackrel{\rightarrow }{%
\alpha }$ is given by 
\[
\stackrel{\rightarrow }{\alpha }(z;X,Y):=\inf \{\alpha (z;E,F)\mid E\in 
\text{FIN}(X)\text{, }F\in \text{FIN}(Y)\text{, }z\in E\otimes F\}\text{,} 
\]
and the {\it {cofinite hull}}\/ $\stackrel{\leftarrow }{\alpha }$ of $\alpha 
$ is given by 
\[
\hspace{0.2cm}\stackrel{\leftarrow }{\alpha }(z;X,Y):=\sup \{\alpha
(Q_K^X\otimes Q_L^Y(z);X/K,Y/L)\mid K\in \text{COFIN}(X)\text{, }L\in \text{
COFIN}(Y)\}\text{.} 
\]
$\alpha $ is called {\it {finitely generated}}\/ if $\alpha $ $=$ $\stackrel{%
\rightarrow }{\alpha }$, {\it {cofinitely generated}}\/ if $\alpha $ $=$ $%
\stackrel{\leftarrow }{\alpha }$ (it is always true that $\stackrel{%
\leftarrow }{\alpha }$ $\leq \alpha $ $\leq $ $\stackrel{\rightarrow }{%
\alpha }$). $\alpha $ is called {\it {right--accessible}} if $\stackrel{%
\leftarrow }{\alpha }$$(z;E,Y)$ $=$ $\stackrel{\rightarrow }{\alpha }(z;E,Y)$
for all $(E,Y)\in $ FIN $\times $ BAN, {\it {left--accessible}} if $%
\stackrel{\leftarrow }{\alpha }$$(z;X,F)$ $=$ $\stackrel{\rightarrow }{%
\alpha }$$(z;X,F)$ for all $(X,F)\in $ BAN $\times $ FIN, and {\it {%
accessible}} if it is right--accessible and left--accessible. $\alpha $ is
called {\it {totally accessible}} if $\stackrel{\leftarrow }{\alpha }$ $=$ $%
\stackrel{\rightarrow }{\alpha }$. The injective norm $\varepsilon $ is
totally accessible, the projective norm $\pi $ is accessible -- but not
totally accessible, and Pisier's construction implies the existence of a
(finitely generated) tensor norm which is neither left-- nor
right--accessible (see \cite{df}, 31.6).

There exists a powerful one--to--one correspondence between finitely
generated tensor norms and maximal Banach ideals which links thinking in
terms of operators with ''tensorial'' thinking and which allows to transfer
notions in the ''tensor language'' to the ''operator language'' and
conversely. We refer the reader to \cite{df} and \cite{oe1} for detailed
informations concerning this subject. Let $X,Y$ be Banach spaces and $%
z=\sum\limits_{i=1}^nx_i^{\prime }\otimes y_i$\ be an Element in $X^{\prime
}\otimes Y$. Then $T_z(x):=\sum\limits_{i=1}^n\langle x,x_i^{\prime }\rangle 
$ $y_i$ defines a finite rank operator $T_z\in {\frak F}(X,Y)$ which is
independent of the representation of $z$ in $X^{\prime }\otimes Y$. Let $%
\alpha $ be a finitely generated tensor norm and $({\frak A},{\bf A})$ be a
maximal Banach ideal. $\alpha $ and $({\frak A},{\bf A})$ are said to be 
{\it {associated}}, notation: 
\[
({\frak A},{\bf A})\sim \alpha \hspace{0.15cm}\text{(shorthand: }{\frak A}%
\sim \alpha \text{, resp.}\hspace{0.1cm}\alpha \sim {\frak A}\text{),} 
\]
if for all $E,F\in $ FIN 
\[
{\frak A}(E,F)=E^{\prime }{\otimes }_\alpha F 
\]
holds isometrically: ${\bf A}(T_z)=\alpha (z;E^{\prime },F)$.

Since we will use them throughout in this paper, let us recall the important
notions of the conjugate operator ideal (cf. \cite{glr}, \cite{jo} and \cite
{oe2}) and the adjoint operator ideal (all details can be found in the
standard references \cite{df} and \cite{p1}). Let $({\frak A},{\bf A})$ be a
quasi--Banach ideal.

\begin{itemize}
\item  Let ${\frak A}^\Delta (X,Y)$ be the set of all $T\in {\cal L}(X,Y)$
which satisfy 
\[
{\bf A}^\Delta (T):=\sup \{\mid tr(TL)\mid \text{ }\mid L\in {\cal F}(Y,X),%
{\bf A}(L)\leq 1\}<\infty .
\]
Then a Banach ideal $({\frak A}^\Delta ,{\bf A}^\Delta )$ is obtained (here, 
$tr(\cdot )$ denotes the usual trace for finite rank operators). It is
called the {\it {conjugate ideal}} of $({\frak A},{\bf A})$.

\item  Let ${\frak A}^{*}(X,Y)$ be the set of all $T\in {\cal L}(X,Y)$ which
satisfy 
\[
{\bf A}^{*}(T):=\sup \{\mid tr(TJ_E^XSQ_K^Y\}\mid \text{ }\mid E\in
FIN(X),K\in COFIN(Y),{\bf A}(S)\leq 1\}<\infty .
\]
Then a Banach ideal $({\frak A}^{*},{\bf A}^{*})$ is obtained. It is called
the {\it adjoint operator ideal }of $({\frak A},{\bf A})$.
\end{itemize}

By definition, it immediately follows that ${\frak A}^\Delta \stackrel{1}{%
\subseteq }{\frak A}^{*}$. Another easy, yet important observation is the
following: let $({\frak A},{\bf A})$ be a quasi--Banach ideal and $({\frak B}%
,{\bf B})$ be a quasi--Banach ideal. If ${\frak A}\stackrel{\frak E}{%
\subseteq }{\frak B}$, then ${\frak B}^{*}\stackrel{1}{\subseteq }{\frak A}%
^{*}$, and ${\frak A}\stackrel{\frak F}{\subseteq }{\frak B}$ implies the
inclusion ${\frak B}^\Delta \stackrel{1}{\subseteq }{\frak A}^\Delta $. In
particular, it follows that ${\frak A}^{\Delta *}\stackrel{1}{=}{\frak A}%
^{**}$ and $({\frak A}^{\Delta \Delta })^{*}\stackrel{1}{=}{\frak A}^{*}$.

In addition to the maximal Banach ideal $({\cal L},{\cal \Vert \cdot \Vert }%
)\sim {\cal \varepsilon }$ we mainly will be concerned with the maximal
Banach ideals $({\frak I},{\bf I})\sim {\cal {\bf \pi }}$ (integral
operators), $({\frak L}_2,{\bf L}_2)\sim w_2$ (Hilbertian operators), $(%
{\frak D}_2,{\bf D}_2)\stackrel{1}{=}({\frak L}_2^{*},{\bf L}_2^{*})%
\stackrel{1}{=}{\frak P}_2^d\circ {\frak P}_2\sim w_2^{*}$\/ ($2$--dominated
operators), $({\frak P}_p,{\bf P}_p)\sim g_p\backslash =g_q^{*}$ (absolutely 
$p$--summing operators), $1\leq p\leq \infty ,\frac 1p+\frac 1q=1$, $({\frak %
L}_\infty ,{\bf L}_\infty )\stackrel{1}{=}({\frak P}_1^{*},{\bf P}%
_1^{*})\sim w_\infty $ and $({\frak L}_1,{\bf L}_1)\stackrel{1}{=}({\frak P}%
_1^{*d},{\bf P}_1^{*d})\sim w_1$. We also consider the maximal Banach ideals 
$({\frak C}_2,{\bf C}_2)$ $\sim c_2$ (cotype 2 operators) and $({\frak A}_P,%
{\bf A}_P)$ $\sim \alpha _P$ (Pisier`s counterexample of a maximal Banach
ideal which is neither right-- nor left--accessible (cf. \cite{df}, 31.6)).

What about the regularity of conjugate ideals? We do not treat this problem
in its whole generality in this paper. In the next section, we will include
additional methods and tools which are of local nature, like accessibility
or the principle of local reflexivity for operator ideals to prove the
regularity of conjugate operator ideals of type ${\frak A}^{\Delta \Delta }$
(cf. Proposition 3.2). However, if $({\frak A},{\bf A})$ is a maximal Banach
ideal, then ${\frak A}^\Delta $ is regular, since:

\begin{proposition}
Let $({\frak A},{\bf A})$ be a quasi--Banach ideal. If ${\frak A}\stackrel{1%
}{=}{\frak A}^{dd}$, then ${\frak A}^\Delta $ is regular.
\end{proposition}

{\sc Proof:} Let $X,Y$ be arbitrary Banach spaces, $T\in {\frak A}^{\Delta
reg}(X,Y)$ and $L\in {\frak F}(Y,X)$. Choose $A\in {\frak F}(Y^{\prime
\prime },X)$ so that $L^{\prime \prime }=j_XA$ (if $L=T_z$ with $%
z=\sum\limits_{i=1}^ny_i^{\prime }\otimes x_i\in Y^{\prime }\otimes X$, then 
$A=T_w$ where $w:=\sum\limits_{i=1}^nj_{Y^{\prime }}y_i^{\prime }\otimes x_i$%
). Since ${\frak A}\stackrel{1}{=}{\frak A}^{dd}$ in particular is regular,
we have 
\begin{eqnarray*}
\left| tr(TL)\right| &=&\left| tr(T^{\prime \prime }j_XA)\right| =\left|
tr(j_YTA)\right| \\
&\leq &{\bf A}^\Delta (j_YT)\cdot {\bf A}(A) \\
&=&{\bf A}^{\Delta reg}(T)\cdot {\bf A}(L)\text{ ,}
\end{eqnarray*}
and the claim follows. $\blacksquare $

Given quasi--Banach ideals $({\frak A},{\bf A)}$ and $({\frak B},{\bf B)}$,
let $({\frak A}\circ {\frak B},{\bf A\circ B})$ be the corresponding product
ideal and $({\frak A}\circ {\frak B}^{-1},{\bf A\circ B}^{-1})$ (resp. $(%
{\frak A}^{-1}\circ {\frak B},{\bf A}^{-1}{\bf \circ B}))$ the corresponding
''right--quotient'' (resp. ''left--quotient''). We write $({\frak A}^{inj},%
{\bf A}^{inj})$, to denote the {\it injective hull}{\em \ }of ${\frak A}$,
the unique smallest injective quasi--Banach ideal which contains $({\frak A},%
{\bf A})$, and $({\frak A}^{sur},{\bf A}^{sur})$, the {\it surjective hull}
of ${\frak A}$, is the unique smallest surjective quasi--Banach ideal which
contains $({\frak A},{\bf A})$. Of particular importance are the quotients $%
{\frak A}^{\dashv }:={\frak I}\circ {\frak A}^{-1}$ and ${\frak A}^{\vdash
}:={\frak A}^{-1}\circ {\frak I}$ and their relations to ${\frak A}^\Delta $
and ${\frak A}^{*}$, treated in detail in \cite{oe1} and \cite{oe4}. Very
useful will be the following statement which represents the injective hull
(resp. the surjective hull) of a conjugate operator ideal as a quotient (cf.
Corollary 3.4):

\begin{proposition}
Let $({\frak A},{\bf A})$ be an arbitrary quasi--Banach ideal. Then 
\[
({\frak A}^\Delta )^{inj}\stackrel{1}{=}{\frak P}_1\circ {\frak A}^{-1}
\]
and 
\[
({\frak A}^\Delta )^{sur}\stackrel{1}{=}{\frak A}^{-1}\circ {\frak P}_1^d
\]
\end{proposition}

{\sc Proof:} It is sufficient to prove the statement only for the injective
hull. Since 
\[
({\frak A}^\Delta )^{inj}\circ {\frak A}\stackrel{1}{\subseteq }({\frak A}%
^\Delta \circ {\frak A})^{inj}\stackrel{1}{\subseteq }{\frak I}^{inj}%
\stackrel{1}{=}{\frak P}_1\text{,} 
\]
{\sc \ }it follows that $({\frak A}^\Delta )^{inj}\stackrel{1}{\subseteq }%
{\frak P}_1\circ {\frak A}^{-1}$. To see the other inclusion, note that 
\[
{\frak A}^\Delta (\cdot ,Y_0)\stackrel{1}{=}{\frak I}\circ {\frak A}%
^{-1}(\cdot ,Y_0) 
\]
holds for every Banach space $Y_0$ of which the dual has the metric
approximation property (this follows by an direct application of \cite{p1},
Lemma 10.2.6.). Hence, 
\[
{\frak P}_1\circ {\frak A}^{-1}\stackrel{1}{=}{\frak I}^{inj}\circ {\frak A}%
^{-1}\stackrel{1}{\subseteq }({\frak I}\circ {\frak A}^{-1})^{inj}\stackrel{1%
}{=}({\frak A}^\Delta )^{inj}\text{,} 
\]
and the proof is finished.$\blacksquare $

A deeper investigation of relations between the Banach ideals $({\frak A}%
^\Delta ,{\bf A}^\Delta )$ and $({\frak A}^{*},{\bf A}^{*})$ needs the help
of an important local property, known as accessibility, which can be viewed
as a local version of injectivity and surjectivity. All necesary details
about accessibility and its applications can be found in \cite{df}, \cite
{oe2}, \cite{oe3} and \cite{oe4}. So let us recall :

\begin{itemize}
\item  A quasi--Banach ideal $({\frak A},{\bf A})$ is called {\it {%
right--accessible}}, if for all $(E,Y)\in $ FIN $\times $ BAN, operators $%
T\in {\cal L}(E,Y)$ and $\varepsilon >0$ there are $F\in $ FIN$(Y)$ and $%
S\in {\cal L}(E,F)$ so that $T=J_F^YS$ and ${\bf A}(S)\leq (1+\varepsilon )%
{\bf A}(T)$.

\item  $({\frak A},{\bf A})$ is called {\it {left--accessible}}, if for all $%
(X,F)\in $ BAN $\times $ FIN, operators $T\in {\cal L}(X,F)$ and $%
\varepsilon >0$ there are $L\in $ COFIN$(X)$ and $S\in {\cal L}(X/L,F)$ so
that $T=SQ_L^X$ and ${\bf A}(S)\leq (1+\varepsilon ){\bf A}(T)$.

\item  A left--accessible and right--accessible quasi--Banach ideal is called%
{\ }{\it {accessible}}.

\item  $({\frak A},{\bf A})$ is {\it {totally accessible}}, if for every
finite rank operator $T\in {\cal F}(X,Y)$ acting between Banach spaces $X$, $%
Y$ and $\varepsilon >0$ there are $(L,F)\in $ COFIN$(X)\times $ FIN$(Y)$ and 
$S\in {\cal L}(X/L,F)$ so that $T=J_F^YSQ_L^X$ and ${\bf A}(S)\leq
(1+\varepsilon ){\bf A}(T)$.
\end{itemize}

Due to the existence of Banach spaces without the approximation property, we
will see now that conjugate hulls are not ''big enough'' to contain such
spaces. To this end, consider an arbitrary {\it Banach} ideal $({\frak A},%
{\bf A})$, and let $X$ be a Banach space so that $Id_X\in {\frak A}^\Delta $
(i.e., $X\in $ space$({\frak A}^\Delta )$). Since $({\frak N},{\bf N})$, the
collection of all nuclear operators, is the smallest Banach ideal, it
follows that $Id_X\in {\frak N}^\Delta $ and ${\bf N}^\Delta (Id_X)\leq {\bf %
A}^\Delta (Id_X)$. Hence, if $T\in {\cal L}(X,X)$ is an arbitrary linear
operator, it follows that $T=TId_X\in {\frak N}^\Delta (X,X)$ and ${\bf N}%
^\Delta (T)\leq \Vert T\Vert \cdot {\bf N}^\Delta (Id_X)\leq \Vert T\Vert
\cdot {\bf A}^\Delta (Id_X)$. But this implies that 
\[
{\cal L}(X,X)={\frak N}^\Delta (X,X).
\]
If ${\frak A}$ contains the class ${\frak I}$ of all integral operators
(e.g., if ${\frak A}$ is maximal or if ${\frak A}$ is a conjugate of a
quasi--Banach ideal), similar considerations lead to 
\[
{\cal L}(X,X)={\frak I}^\Delta (X,X)\text{,}
\]
and \cite[Proposition 2.2.]{jo} now imply the following

\begin{remark}
Let $({\frak A},{\bf A})$ be an arbitrary quasi--Banach ideal, and let $X$
be a Banach space so that $X\in $ space$({\frak A}^\Delta )$. If ${\frak A}$
is normed, then $X$ has the approximation property. If ${\frak I}\subseteq 
{\frak A}$, then $X$ has the bounded approximation property.
\end{remark}

\begin{corollary}
Let $({\frak A},{\bf A})$ be an arbitrary quasi--Banach ideal so that there
exists a Banach space in space$({\frak A}^{**})$ without the bounded
approximation property, then ${\frak A}^\Delta $ (-- in particular ${\frak A}%
^{*}$) cannot be totally accessible.
\end{corollary}

{\sc Proof:} Let $X$ be a Banach space without the bounded approximation
property so that $X\in $ space$({\frak A}^{**})$. Assume, ${\frak A}^\Delta $
is totally accessible, then 
\[
{\frak A}^{**}\stackrel{1}{=}{\frak A}^{\Delta *}\stackrel{1}{=}{\frak A}%
^{\Delta \Delta } 
\]
Since ${\frak I}\stackrel{1}{=}{\frak L}^\Delta \stackrel{1}{\subseteq }%
{\frak A}^\Delta $, the previous Remark leads to a contradiction\footnote{%
Proposition 21.6 in \cite{df} is a special case of this Corollary.}.$%
\blacksquare $

Since {\it ultrastable }operator ideals play an important part in this
paper, we completely recall the definition of an ultrastable operator ideal
and its construction (cf. \cite{dk}, \cite{df}, \cite{djt}, \cite{h}, \cite
{k} and \cite{p1}): Let $I$ be a non--empty set and ${\cal U}$ be an
ultrafilter in $I$. If $(X_i)_{i\in I}$ is a family of Banach spaces,
consider in the Banach space 
\[
l_\infty (X_i;I):=\left\{ x=(x_i)_{i\in I}\in \prod_{i\in I}X_i\;\mid
\;\Vert x\Vert _\infty :=\sup_{i\in I}\Vert x_i\Vert <\infty \right\} 
\]
the closed subspace $N_{{\cal U}}(X_i;I):=\left\{ (x_i)_{i\in I}\;\in
l_\infty (X_i;I)\;\mid \;\lim_{{\cal U}}\Vert x_i\Vert =0\;\right\} $. The 
{\it ultraproduct }of the family $(X_i)_{i\in I}$ with respect to the
ultrafilter{\it \ }${\cal U}$ is defined to be the Banach space 
\[
(\prod_{i\in I}X_i)_{{\cal U}}:=l_\infty (X_i;I)\diagup N_{{\cal U}}(X_i;I) 
\]
equipped with the canonical quotient norm. The elements of $(\prod_{i\in
I}X_i)_{{\cal U}}$ are denoted by $(x_i)_{{\cal U}}$ (whenever $(x_i)_{i\in
I}\in l_\infty (X_i;I)$), and the construction implies that $\Vert (x_i)_{%
{\cal U}}\Vert =\;\lim_{{\cal U}}\Vert x_i\Vert $. If $(X_i)_{i\in I}$ and $%
(Y_i)_{i\in I}$ are two families of Banach spaces and $T_i\in {\frak L}%
(X_i,Y_i)$ with $c:=\sup_{i\in I}\Vert T_i\Vert <\infty $, then 
\[
T^I(x_i):=(T_ix_i) 
\]
defines an operator $T^I:l_\infty (X_i;I)\longrightarrow l_\infty (Y_i;I)$
with $\Vert T^I\Vert \leq c$ and which maps $N_{{\cal U}}(X_i;I)$ into $N_{%
{\cal U}}(Y_i;I)$. Consequently, there exists a linear operator $%
T:(\prod_{i\in I}X_i)_{{\cal U}}\longrightarrow (\prod_{i\in I}Y_i)_{{\cal U}%
}$ so that $T(x_i)_{{\cal U}}=(T_ix_i)_{{\cal U}}$. $T$ is called the {\it %
ultraproduct }$(T_i)_{{\cal U}}$ {\it of the operators} $T_i$. It satisfies $%
\Vert (T_i)_{{\cal U}}\Vert =\lim_{{\cal U}}\Vert T_i\Vert $.

Let $({\frak A},{\bf A})$ be an arbitrary quasi--Banach ideal. ${\frak A}$
is called {\it ultrastable} if for every ultrafilter ${\cal U}$ on $I$ and
every ${\bf A}$--bounded family of operators $T_i\in {\frak A}(X_i,Y_i)$ 
\[
(T_i)_{{\cal U}}\in {\frak A}((\prod_{i\in I}X_i)_{{\cal U}},(\prod_{i\in
I}Y_i)_{{\cal U}})\text{ \ and \ }{\bf A(}(T_i)_{{\cal U}})\leq \lim_{{\cal U%
}}{\bf A}(T_i)\text{ .} 
\]
The key part of ultrastable operator ideals is given by the following
relation (see \cite{p1}, Theorem 8.8.6.):

\begin{theorem}[Pietsch]
Let $({\frak A},{\bf A})$ be an ultrastable quasi--Banach ideal. Then 
\[
({\frak A},{\bf A})^{\max }=({\frak A},{\bf A})^{reg}\text{ .}
\]
\end{theorem}

Although ${\frak A}^{min}$ (resp. $({\frak A}^{*\Delta })^{dd}$) is always
accessible, Pisier's counterexample shows the existence of maximal Banach
ideals which neither are left nor right--accessible. However, accessibility
conditions of a quasi--Banach ideal at least can be transmitted to its
regular hull:

\begin{proposition}
Let $({\frak A},{\bf A)}$ be an arbitrary quasi--Banach ideal. If ${\frak A}$
is right--accessible (resp. totally--accessible), then the regular hull $%
{\frak A}^{reg}$ is also right--accessible (resp. totally--accessible).
\end{proposition}

{\sc Proof: }Let $\epsilon >0$, $X$, $Y$ be Banach spaces and $%
T\in {\frak F}(X,Y)$ an arbitrary finite rank operator. Assume that ${\frak A%
}$ is totally accessible or that $X\in FIN$ and ${\frak A}$ is
right--accessible. In both cases, there exists a finite dimensional Banach
space $F\in FIN(Y^{\prime \prime })$ and an operator $S\in {\frak L}(X,F)$,
so that $j_YT=J_F^{Y^{\prime \prime }}S$ and 
\[
{\bf A}(S)<(1+\epsilon )\cdot {\bf A}(j_YT)=(1+\epsilon )\cdot {\bf A}%
^{reg}(T)\text{.} 
\]
Due to the classical principle of local reflexivity for linear operators,
there exists an operator $W\in {\frak L}(F,Y)$ so that ${\bf \Vert }W{\bf %
\Vert }<1+\epsilon $ and $j_YWz=J_F^{Y^{\prime \prime }}z$ for all $z\in F$
which satisfy $J_F^{Y^{\prime \prime }}z\in j_Y(Y)$. Let $x\in X$ and put $%
z:=Sx$. Then $J_F^{Y^{\prime \prime }}z=$ $j_YTx\in j_Y(Y)$, which therefore
implies that $j_YWSx=J_F^{Y^{\prime \prime }}z=$ $j_YTx$. Now, factor $W$
canonically through a finite dimensional subspace $G$ of $Y$ so that $%
W=J_G^YU$ and ${\bf \Vert }U{\bf \Vert }<1+\epsilon $. Consequently, $%
T=WS=J_G^Y(US)$, and 
\[
{\bf A}^{reg}(US)<(1+\epsilon )^2\cdot {\bf A}^{reg}(T)\text{.} 
\]
Hence, ${\frak A}^{reg}$ is right--accessible (in each of the both cases).
In the case of ${\frak A}$ being totally accessible, the operator $S$ even
can be chosen as $S=S_0Q_K^X$, where $K\in COFIN(X)$ and $S_0\in {\frak L}%
(X\diagup K,F)$ so that 
\[
{\bf A(}S_0)<(1+\epsilon )\cdot {\bf A}^{reg}(T)\text{,} 
\]
and the proof is finished.$\blacksquare $

Let us finish this section with a short Remark concerning ultrastability
versus accessibility. To this end, let $({\frak A},{\bf A})$ be a maximal
Banach ideal. Then ${\frak A}$ is right--accessible (resp. totally
accessible) {\it if and only if} ${\frak A}^{*}\stackrel{1}{=}{\frak I}\circ 
{\frak A}^{-1}$ (resp. ${\frak A}^{*}\stackrel{1}{=}{\frak A}^\Delta $) (cf. 
\cite{oe4}). A straightforward calculation therefore leads to the following
structurally interesting

\begin{remark}
Let $({\frak A},{\bf A})$ be a maximal Banach ideal. Then the following
statements are equivalent:

\begin{enumerate}
\item[(i)]  ${\frak A}$ is right--accessible (resp. totally accessible)

\item[(ii)]  ${\frak I}\circ {\frak A}^{-1}$(resp. ${\frak A}^\Delta $) is
ultrastable.
\end{enumerate}
\end{remark}

\section{Extension of finite rank operators and the principle of local
reflexivity for operator ideals}

Let $({\frak A},{\bf A})$ be a {\it maximal} Banach ideal. Then, ${\frak A}%
^\Delta $ always is right--accessible (cf. \cite{oe4}). The natural question
whether ${\frak A}^\Delta $ is {\it {left}}--accessible is still open$%
\footnote{%
For minimal Banach ideals $({\frak A},{\bf A)}$, there exist
counterexamples: The conjugate of ${\frak A}_P^{\min }$ neither is
right--accessible nor left--accessible (cf. \cite{oe2}).}$ and leads to
interesting and non--trivial results concerning the local structure of $%
{\frak A}^\Delta $. Deeper investigations of the left--accessibility of $%
{\frak A}^\Delta $ namely lead to a link with a principle of local
reflexivity for operator ideals (a detailed discussion can be found in \cite
{oe1} and \cite{oe2}) which allows a transmission of the operator norm
estimation in the classical principle of local reflexivity to the ideal norm 
${\bf A}$. So let us recall the

\begin{definition}
Let $E$ and $Y$ be Banach spaces, $E$ finite dimensional, $F\in $ FIN$%
(Y^{\prime })$ and $T\in {\frak L}(E,Y^{\prime \prime })$. Let $({\frak A},%
{\bf A})$ be a quasi--Banach ideal and $\epsilon >0$. We say that the
principle of ${\frak A}-$local reflexivity (short: ${\frak A}-LRP$) is
satisfied, if there exists an operator $S\in {\frak L}(E,Y)$ so that

\begin{enumerate}
\item[(1)]  ${\bf A}(S)\leq (1+\epsilon )\cdot {\bf A}^{**}(T)$

\item[(2)]  $\left\langle Sx,y^{\prime }\right\rangle =\left\langle
y^{\prime },Tx\right\rangle $ for all $(x,y^{\prime })\in E\times F$

\item[(3)]  $j_YSx=Tx$ for all $x\in T^{-1}(j_Y(Y))$.
\end{enumerate}
\end{definition}

Although both, the quasi--Banach ideal ${\frak A}$ and{\it \ }the $1$%
--Banach ideal ${\frak A}^{**}$ are involved, the asymmetry can be justified
by the following statement which holds for arbitrary quasi--Banach ideals
(see \cite{oe2}):

\begin{theorem}
Let $({\frak A},{\bf A})$ be a quasi--Banach ideal. Then the following
statements are equivalent:

\begin{enumerate}
\item[(i)]  ${\frak A}^\Delta $ is left--accessible

\item[(ii)]  ${\frak A}^{**}(E,Y^{\prime \prime })$ $\widetilde{=}$ ${\frak A%
}(E,Y)^{\prime \prime }$ for all $(E,Y)\in $ FIN $\times $ BAN

\item[(iii)]  The ${\frak A}-LRP$ holds.
\end{enumerate}
\end{theorem}

One reason which leads to extreme persistent difficulties concerning the
verification of the ${\frak A}-LRP$ for a given maximal Banach ideal ${\frak %
A}$, is the behaviour of the bidual $({\frak A}^\Delta )^{dd}$: although we
know that in general $({\frak A}^\Delta )^{dd}$ is accessible (see \cite{oe1}
and \cite{oe2}) and that $({\frak A}^\Delta )^{dd}\stackrel{1}{\subseteq }%
{\frak A}^\Delta $, we do not know whether ${\frak A}^\Delta (X,Y)$ and $(%
{\frak A}^\Delta )^{dd}(X,Y)$ coincide {\it isometrically} for {\it all}
Banach spaces $X$ and $Y$. If we allow in addition the approximation
property of $X$ or $Y$, then we may state the following

\begin{lemma}
Let $({\frak A},{\bf A})$ be an arbitrary maximal Banach ideal and $X$, $Y$
be arbitrary Banach spaces. Then 
\[
{\frak A}^{d\Delta }(X,Y)\stackrel{1}{=}{\frak A}^{\Delta d}(X,Y)
\]
holds in each of the following two cases:

\begin{enumerate}
\item[(i)]  $X^{\prime }$ has the metric approximation property

\item[(ii)]  $Y^{\prime }$ has the metric approximation property and the $%
{\frak A}^d-LRP$ is satisfied.
\end{enumerate}
\end{lemma}

{\sc Proof:} Only the inclusion $\subseteq $ is not trivial. So, let $T\in 
{\frak A}^{d\Delta }(X,Y)$ be given. First, we consider the case (i). Due to
Proposition 2.3 of \cite{jo}, it follows that in general 
\[
{\frak A}^{d\Delta }\stackrel{1}{\subseteq }({\frak A}^d)^{-1}\circ {\frak I}%
\stackrel{1}{=}({\frak A}^{\dashv })^d\text{,} 
\]
so that $T^{\prime }\in {\frak A}^{\dashv }(Y^{\prime },X^{\prime })$, and $%
{\bf A}^{\dashv }(T^{\prime })\leq {\bf A}^{d\Delta }(T)$. Since $X^{\prime
} $ has the metric approximation property we even obtain that $T^{\prime
}=Id_{X^{\prime }}T^{\prime }\in {\frak I}^\Delta \circ {\frak A}^{\dashv
}(Y^{\prime },X^{\prime })\stackrel{1}{\subseteq }{\frak A}^\Delta
(Y^{\prime },X^{\prime })$, and case (i) is finished.

To prove case (ii), we have to proceed in a total different way.
Let $L\in {\frak F}(X^{\prime },Y^{\prime })$ be an arbitrary finite rank
operator and $\epsilon >0$. Since $Y^{\prime }$ has the metric approximation
property, there exists a finite rank operator $A\in {\frak F}(Y^{\prime
},Y^{\prime })$ so that $L=AL$ and $\Vert A\Vert \leq 1+\epsilon $. Thanks
to canonical factorization, we can find a finite dimensional space $G$ and
operators $A_1\in {\frak L}(Y^{\prime },G^{\prime \prime })$, $A_2\in {\frak %
L}(G^{\prime \prime },Y^{\prime })$ so that $A=A_2A_1$, $\Vert A_2\Vert \leq
1$ and $\Vert A_1\Vert \leq 1+\epsilon $. Now, look carefully at the
composition of the two operators $A_1L\in {\frak F}(X^{\prime },G^{\prime
\prime })$ and $T^{\prime }A_2\in {\frak L}(G^{\prime \prime },X^{\prime })$%
. Using exactly the same considerations as in \cite[E.3.2.]{p1}, the assumed 
${\frak A}^d-LRP$ implies the existence of an operator $\Lambda \in {\frak L}%
(G^{\prime },X)$ so that$\footnote{%
We only have to substitute the operator norm through the ideal norm ${\bf A}%
^d$.}$%
\[
{\bf A}^d(\Lambda )\leq (1+\epsilon )\cdot {\bf A}^d((A_1L)^{\prime
})=(1+\epsilon )\cdot {\bf A}(A_1L)
\]
and $A_1LT^{\prime }A_2=\Lambda ^{\prime }T^{\prime }A_2$. Since $G$ is
finite dimensional, we may represent $A_2$ as the dual of a finite rank
operator $B_2\in {\frak F}(Y,G^{\prime })$, and consequently it follows 
\begin{eqnarray*}
&\mid &tr(T^{\prime }L)\mid =\mid tr(A_1LT^{\prime }A_2)\mid =\mid
tr(\Lambda ^{\prime }T^{\prime }A_2)\mid =\mid tr(T\Lambda B_2)\mid  \\
&\leq &{\bf A}^{d\Delta }(T)\cdot {\bf A}^d(\Lambda ) \\
&\leq &(1+\epsilon )^2\cdot {\bf A}^{d\Delta }(T)\cdot {\bf A}(L)\text{.}
\end{eqnarray*}
Hence, $T^{\prime }\in {\frak A}^\Delta (Y^{\prime },X^{\prime })$, and $%
{\bf A}^\Delta (T^{\prime })\leq {\bf A}^{d\Delta }(T)$, and case (ii) also
is proved.$\blacksquare $

A straightforward dualization of the previous Lemma implies a result which
we will use later again:

\begin{corollary}
Let $({\frak A},{\bf A})$ be an arbitrary maximal Banach ideal and $X,Y$ be
arbitrary Banach spaces. Then 
\[
{\frak A}^\Delta (X,Y)\stackrel{1}{=}({\frak A}^\Delta )^{dd}(X,Y)
\]
holds in each of the following two cases:

\begin{enumerate}
\item[(i)]  $X^{\prime \prime }$ has the metric approximation property and
the ${\frak A}^d-LRP$ is satisfied

\item[(ii)]  $Y^{\prime \prime }$ has the metric approximation property and
the ${\frak A}-LRP$ is satisfied.
\end{enumerate}
\end{corollary}

Using the considerations of the previous section, we obtain a closer
approach to the ${\frak A}-LRP$ in the following sense:

\begin{theorem}
Let $({\frak A},{\bf A)}$ be an arbitrary Banach ideal. If ${\frak A}$ is
right--accessible and ultrastable, then the ${\frak A}-LRP$ is satisfied.
\end{theorem}

{\sc Proof:} Let $\epsilon >0$. Let $E$ and $Y$ be Banach spaces, $%
E$ finite dimensional, $F\in $ FIN$(Y^{\prime })$ and $T\in {\frak L}%
(E,Y^{\prime \prime })$. Since the right--accessibility of ${\frak A}$
implies the right--accessibility of ${\frak A}^{reg}$, there exists $G\in $
FIN$(Y^{\prime \prime })$ and an operator $B\in {\frak L}(E,G)$ so that{\bf %
\ }${\bf A}(B)={\bf A}^{reg}(B)\leq (1+\epsilon )\cdot {\bf A}^{reg}(T)$ and 
$T=J_G^{Y^{\prime \prime }}B$. The (classical) principle of local
reflexivity, applied to the operator $J_G^{Y^{\prime \prime }}$, implies the
existence of a further operator $\Lambda \in {\frak L}(G,Y)$ so that $\Vert
\Lambda \Vert \leq 1+\epsilon $ and $\left\langle \Lambda z,y^{\prime
}\right\rangle =\left\langle y^{\prime },J_G^{Y^{\prime \prime
}}z\right\rangle $ for all $(z,y^{\prime })\in G\times F$. Hence, $%
S:=\Lambda B\in {\frak L}(E,Y)$, 
\begin{equation}
{\bf A}(S)\leq (1+\epsilon )^2\cdot {\bf A}^{reg}(T)\text{,}  \tag{$*$}
\end{equation}
and $\left\langle Sx,y^{\prime }\right\rangle =\left\langle y^{\prime
},J_G^{Y^{\prime \prime }}(Bx)\right\rangle =\left\langle y^{\prime
},Tx\right\rangle $ for all $(x,y^{\prime })\in E\times F$. Now, the
assumption further implies that ${\bf A}^{reg}(T)={\bf A}^{\max }(T)={\bf A}%
^{**}(T)$, and the proof is finished.$\blacksquare $

Consequently, every right--accessible and maximal Banach ideal $({\frak A},%
{\bf A})$ satisfies the ${\frak A}-LRP$. Is the converse implication also
true? Does the ${\frak A}^{**}-LRP$ even imply the right--accessibility of $%
{\frak A}^{**}$? A partial answer -- involving Hilbert spaces -- is given in
Corollary 4.2. For minimal operator ideals, the previous considerations
immediately lead to the following fact which we want to state separately:

\begin{corollary}
Let $({\frak A},{\bf A)}$ be an arbitrary Banach ideal. If ${\frak A}^{\min }
$ is ultrastable, then the ${\frak A}-LRP$ is satisfied.
\end{corollary}

{\sc Proof:} Since ${\frak A}$ is a Banach ideal, ${\frak A}^{\min }$ is
also a Banach ideal (cf. \cite{b}, 9.2. and \cite{df}, 22.2.), and ${\frak A}%
^{\min }$ always is (right--)accessible. Therefore, the assumed
ultrastability of ${\frak A}^{\min }$, implies the validity of the ${\frak A}%
^{\min }-LRP$ and in particular the validity of the ${\frak A}-LRP$.$%
\blacksquare $

Although operator ideals which are both, minimal {\it and} ultrastable, seem
to be quite strange objects, there exist examples, such as the quasi--Banach
ideal $({\frak N}_{(r,p,q)},{\bf N}_{(r,p,q)})$ (the collection of all $%
(r,p,q)$--nuclear operators). If $0<r<\infty $, $1\leq p,q\leq \infty $ and $%
1+1/r>1/p+1/q$, then ${\frak N}_{(r,p,q)}$ is ultrastable and minimal (see 
\cite{p1}, 18.1.4. and 18.1.9.).

Pisier`s counterexample of the maximal Banach ideal $({\frak A}_P,{\bf A}_P)$
which neither is left--accessible nor right--accessible (cf. \cite{df},
31.6) implies that in particular $({\frak A}_P^{*},{\bf A}_P^{*})$ neither
is left--accessible nor right--accessible. Thinking at ${\frak A}_P^\Delta 
\stackrel{1}{\subseteq }{\frak A}_P^{*}$, this leads to the natural and even
more tough question whether the ${\frak A}_P-LRP$ is true or false. However,
due to Corollary 2.1, we already know that ${\frak A}_P^\Delta $ cannot be
totally{\it \ }accessible. Is it even true that $({\frak A}_P^\Delta )^{inj}%
\stackrel{1}{=}{\frak P}_1\circ ({\frak A}_P)^{-1}$ is not totally
accessible? If this is the case, the ${\frak A}_P-LRP$ will be false.
Unfortunately, we will later recognize that, in addition, ${\frak A}%
_P^\Delta $ (and ${\frak I}\circ ({\frak A}_P)^{-1}$) cannot be injective.
What about the left accessibility of ${\frak A}_P^{*\Delta }$? Although we
do not investigate the local structure of ${\frak A}_P$ in this paper, we
want to show a way how to construct other counterexamples. A first step
towards an construction of such a candidate $({\frak A},{\bf A})$ is given
by the following factorization property for finite rank operators which had
been introduced by Jarchow and Ott in their paper \cite{jo}. It not only
turns out to be a useful tool in constructing such a counterexample; later,
we will also use this factorization property to show that ${\frak L}_\infty $
{\it is not totally accessible} -- answering an open question of Defant and
Floret (see \cite{df}, 21.12)! So let us recall the definition of this
factorization property and its implications:

\begin{definition}[Jarchow/Ott]
Let $({\frak A},{\bf A)}$ and $({\frak B},{\bf B)}$ be arbitrary
quasi--Banach ideals. Let $L\in {\frak F}(X,Y)$ an arbitrary finite rank
operator between two Banach spaces $X$ and $Y$. Given $\epsilon >0$, we can
find a Banach space $Z$ and operators $A\in {\frak A}(Z,Y)$, $B\in {\frak B}%
(X,Z)$ so that $L=AB$ and 
\[
{\bf A}(A)\cdot {\bf B}(B)\leq (1+\epsilon )\cdot {\bf A\circ B(}L)\text{.}
\]

\begin{enumerate}
\item[(i)]  If the operator $A$ is of finite rank, we say that ${\frak A}%
\circ {\frak B}$ has the property{\it \ }(I).

\item[(ii)]  If the operator $B$ is of finite rank, we say that ${\frak A}%
\circ {\frak B}$ has the property{\it \ }(S).
\end{enumerate}
\end{definition}

Important examples are the following (see \cite{jo}, Lemma 2.4.):

\begin{itemize}
\item  If ${\frak B}$ is injective, or if ${\frak A}$ contains ${\frak L}_2$
as a factor, then ${\frak A}\circ {\frak B}$ has the property{\it \ }(I).

\item  If ${\frak A}$ is surjective, or if ${\frak B}$ contains ${\frak L}_2$
as a factor, then ${\frak A}\circ {\frak B}$ has the property{\it \ }(S).
\end{itemize}

Since ${\frak L}_2\circ {\frak A}$ is injective for every quasi--Banach
ideal $({\frak A},{\bf A})$ (see \cite{oe4}, Lemma 5.1.), ${\frak B}\circ 
{\frak L}_2\circ {\frak A}$ therefore has the property (I) as well as the
property (S), for all quasi--Banach ideals $({\frak B},{\bf B})$. Such
ideals are exactly those which contain ${\frak L}_2$ as factor -- in the
sense of \cite{jo}.

The next statement will be also useful for our further investigatons (see 
\cite{jo}, 2.5.):

\begin{proposition}
Let $({\frak A},{\bf A)}$ and $({\frak B},{\bf B)}$ be arbitrary
quasi--Banach ideals. Then

\begin{enumerate}
\item[(i)]  $({\frak A}\circ {\frak B)}^\Delta \stackrel{1}{=}{\frak B}%
^{-1}\circ {\frak A}^\Delta $, if ${\frak A}\circ {\frak B}$ has the property%
{\it \ }(I).

\item[(ii)]  $({\frak A}\circ {\frak B)}^\Delta \stackrel{1}{=}{\frak B}%
^\Delta \circ {\frak A}^{-1}$, if ${\frak A}\circ {\frak B}$ has the property%
{\it \ }(S).

In both cases (i) and (ii), the inclusion $\stackrel{1}{\subseteq }$ holds
in general -- without any assumption on the ideals ${\frak A}$ and ${\frak B}
$.
\end{enumerate}
\end{proposition}

Later, we will recognize the particular importance of operator ideals of
type ${\frak A}^{*}\circ {\frak L}_\infty $ which in addition have the
property (I). First, let us note an implication of this factorization
property which gives us a further insight into the local structure of
conjugate operator ideals:

\begin{corollary}
Let $({\frak A},{\bf A)}$ be a quasi--Banach ideal so that ${\frak A}^\Delta 
$ is injective or surjective, then space$({\frak A)}$ cannot contain a
Banach space without the approximation property.
\end{corollary}

{\sc Proof:} Again, a proof for the injective case is enough. So, assume
that the statement is false. Choose a Banach space $X\in $ space$({\frak A)}$
without the approximation property. Since ${\frak A}^\Delta $ is injective,
it follows that ${\frak L}\circ {\frak A}^\Delta $ has the property (I), so
that 
\[
Id_X\in {\frak A}\stackrel{1}{\subseteq }({\frak A}^\Delta )^{-1}\circ 
{\frak I}\stackrel{1}{=}({\frak L}\circ {\frak A}^\Delta )^\Delta \stackrel{1%
}{=}{\frak A}^{\Delta \Delta }\stackrel{1}{\subseteq }{\frak N}^\Delta \text{%
,} 
\]
which is a contradiction.$\blacksquare $

Next, we will see how the property (I) influences the structure of operator
ideals of type ${\frak A}^{inj*}\stackrel{1}{=}\diagdown {\frak A}^{*}$ and
their conjugates. To this end, first note that for all Banach spaces $X$, $Y$
and $X\stackrel{1}{\hookrightarrow }Z$, every operator $T\in ({\frak A}%
^{inj})^{*}(X,Y)\stackrel{1}{=}\diagdown {\frak A}^{*}(X,Y)$ satisfies the
following extension property: Given $\epsilon >0$, there exists an operator $%
\widetilde{T}\in \diagdown {\frak A}^{*}(Z,Y^{\prime \prime })$ so that $%
j_YT=\widetilde{T}J_X^Z$ and $\diagdown {\bf A}^{*}(\widetilde{T})\leq
(1+\epsilon )\cdot \diagdown {\bf A}^{*}(T)$ (see \cite{h1}, Satz 7.14). In
particular, such an extension holds for all finite rank operators. However,
we then cannot be sure that $\widetilde{T}$ is also as a{\it \ finite rank }%
operator. Here, property (I) comes into play -- in the following sense:

\begin{theorem}
Let $({\frak A},{\bf A)}$ be a maximal Banach ideal so that ${\frak A}%
^{*}\circ {\frak L}_\infty $ has the property{\it \ }(I). Let $\epsilon >0$, 
$X$ and $Y$ be arbitrary Banach spaces and $L\in {\frak F}(Y,X)$. Let $Z$ be
a Banach space which contains $Y$ isometrically. Then there exists a finite
rank operator $V\in {\frak F}(Z,X^{\prime \prime })$ so that $j_XL=VJ_Y^Z$
and 
\[
{\bf A}^{*}(V)\leq (1+\epsilon )\cdot ({\bf A}^{inj})^{*}(L)\text{.}
\]
If in addition, the ${\frak A}^{*}-LRP$ is satisfied, then $V$ even can be
chosen to be a finite rank operator with range in $X$ and $L=VJ_Y^Z$.
\end{theorem}

{\sc Proof:} Let $L\in {\frak F}(Y,X)$ be an arbitrary finite rank
operator between arbitrarily given Banach spaces $X$ and $Y$, and set $(%
{\frak B},{\bf B)}:=({\frak A}^{inj},{\bf A}^{inj})$. Let $\epsilon >0$.
Since ${\frak B}^{*}\stackrel{1}{=}({\frak A}^{*}\circ {\frak L}_\infty
)^{reg}$ (cf. \cite{oe4}), there exist a Banach space $W$ and operators $%
A\in {\frak A}^{*}(W,X^{\prime \prime })$, $B\in {\frak L}_\infty (Y,W)$ so
that $j_XL=AB$ and 
\[
{\bf A}^{*}(A)\cdot {\bf L}_\infty (B)\leq (1+\epsilon )\cdot {\bf B}^{*}(L)%
\text{.} 
\]
Due to the assumed property (I) of ${\frak A}^{*}\circ {\frak L}_\infty $,
we even may assume that $A$ is a finite rank operator. Further, we also may
choose a Borel--Radon measure $\mu $ and operators $R\in {\frak L}(L_\infty
(\mu ),W^{\prime \prime })$, $S\in {\frak L}(Y,L_\infty (\mu ))$ so that $%
j_WB=RS$ and 
\[
\Vert R\Vert \cdot \Vert S\Vert \leq (1+\epsilon )\cdot {\bf L}_\infty (B) 
\]
(cf. \cite{df}, 20.12). Due to the metric extension property of $L_\infty
(\mu )$, the operator $S$ can be extended to an operator $\widetilde{S}\in 
{\frak L}(Z,L_\infty (\mu ))$ so that $S=\widetilde{S}J_Y^Z$ and $\Vert 
\widetilde{S}\Vert =\Vert S\Vert $. If we also take into account that $%
Id_{X^{\prime \prime }}^{}=j_{X^{\prime }}^{\prime }j_{X^{\prime \prime
}}^{} $, then we obtain the following factorization of $j_XL$: 
\[
j_XL=j_{X^{\prime }}^{\prime }(j_{X^{\prime \prime }}j_XL)=j_{X^{\prime
}}^{\prime }(A^{\prime \prime }R\widetilde{S}J_Y^Z)\text{.} 
\]
Therefore, $V:=j_{X^{\prime }}^{\prime }A^{\prime \prime }R\widetilde{S}\in 
{\frak F}(Z,X^{\prime \prime })$ is the desired finite rank operator, and
the factorization further shows that 
\[
{\bf A}^{*}(V)\leq (1+\epsilon )^2\cdot {\bf B}^{*}(L)\text{,} 
\]
and the first part of our Theorem is proven.

Now let us assume that in addition the ${\frak A}^{*}-LRP$ is satisfied.
Since $Y$ embeds isometrically into $Y^\infty =l_\infty (B_{Y^{\prime }})$,
the previous considerations (in particular) imply the existence of a finite
rank operator $V\in {\frak F}(Y^\infty ,X^{\prime \prime })$, so that $%
j_XL=VJ_Y$ and ${\bf A}^{*}(V)\leq (1+\epsilon )\cdot {\bf B}^{*}(L)$. Due
to the metric approximation property of the dual of $Y^\infty $, we can find
a finite dimensional subspace $F$ in $Y^\infty $ and an operator $B\in 
{\frak L}(Y^\infty ,F)$ so that $\Vert B\Vert \leq 1+\epsilon $ and $V=WB$
where $W:=VJ_F^{Y^\infty }\in {\frak L}(F,X^{\prime \prime })$. Due to the
assumed ${\frak A}^{*}-LRP$, we even can find an operator $W_0\in {\frak L}%
(F,X)$ so that 
\[
{\bf A}^{*}(W_0)\leq (1+\epsilon )\cdot {\bf A}^{*}(W)\leq (1+\epsilon
)^2\cdot {\bf B}^{*}(L) 
\]
and 
\[
Wx=j_XW_0x\text{ \ for all }x\in W^{-1}(j_X(X))\text{.} 
\]
Since for every $y\in Y$, $x=BJ_Yy\in F$ and $Wx=WBJ_Yy=VJ_Yy=$ $j_XLy\in
j_X(X)$, it therefore follows that 
\[
j_XLy=j_XW_0BJ_Yy\text{\ \ for all }y\in Y\text{ .} 
\]
Hence, $L=W_0BJ_Y$ and ${\bf A}^{*}(W_0B)\leq (1+\epsilon )^3\cdot ({\bf A}%
^{inj})^{*}(L)$. Since $Y^\infty $ has the metric extension property, we can
factorize $J_Y$ as $J_Y=\widetilde{J}J_Y^Z$ so that $\widetilde{J}\in {\frak %
L}(Z,Y^\infty )$, $\Vert \widetilde{J}\Vert =1$, and $V_0:=W_0B\widetilde{J}%
\in {\frak L}(Z,X)$ is our desired finite rank operator.$\blacksquare $

Let $({\frak A},{\bf A)}$ be a Banach ideal and $({\frak A}^{inj},{\bf A}%
^{inj}{\bf )}$ its injective hull. Thinking carefully about the previous
statement, one might guess a strong relationship between the conjugate of $(%
{\frak A}^{inj})^{*}$ and the injective hull of ${\frak A}^{*\Delta }$ --
involving the ${\frak A}^{*}-LRP$ and further accessibility conditions.
Indeed, we will show that such interesting relations exist and that they
even support the search for a counterexample of a maximal Banach ideal $%
{\frak A}_0$ which does not satisfy the ${\frak A}_0-LRP$. So, let us start
with a deeper investigation of the Banach ideal ${\frak A}^{inj*\Delta }$.

\begin{proposition}
Let $({\frak A,}{\bf A})$ be a $1$--Banach ideal so that the ${\frak A}%
^{*}-LRP$ is valid. Then 
\begin{equation}
{\frak A}^{*\Delta inj}\stackrel{\frak F}{=}{\frak A}^{\min inj}\text{ }%
\stackrel{\frak F}{=}{\frak A}^{inj\min }\stackrel{\frak F}{=}{\frak A}%
^{inj*\Delta }\text{ .}  \tag{$*$}
\end{equation}
In particular, ${\frak A}^{inj*\Delta }$ is totally accessible and $({\frak A%
}^{inj*})^{\Delta \Delta }$ regular.
\end{proposition}

{\sc Proof:} Let $\epsilon >0$, $X$ and $Y$ be arbitrary Banach spaces and $%
T\in {\frak F}(X,Y)$ an arbitrary {\it finite rank} operator. Due to the
assumed ${\frak A}^{*}-LRP$, ${\frak A}^{*\Delta }$ is left--accessible, and
it follows the total accessibility of its injective hull $({\frak B},{\bf B}%
):=({\frak A}^{*\Delta inj},{\bf A}^{*\Delta inj})$. Hence, there exist $%
F\in FIN(Y)$, $K\in COFIN(X)$ and an elementary operator $S\in {\frak L}%
(X\diagup K,F)$ so that $T=J_F^YSQ_K^X$ and ${\bf B(}S)<(1+\epsilon )\cdot 
{\bf B}(T)$. Since ${\bf A}$ is an ideal--{\it norm}, we obtain (cf. \cite
{p1}, 8.7.13., 9.2.2. and 9.3.1.) 
\[
{\bf A}^{inj}(S)={\bf A}^{inj**}(S)={\bf A}^{**inj}(S)={\bf B}^{**}(S)={\bf %
B(}S)\text{,} 
\]
so that 
\begin{eqnarray*}
{\bf A}^{inj*\Delta }(T) &\leq &{\bf A}^{inj\min }(T)\leq {\bf A}%
^{inj}(S)\leq {\bf B(}S) \\
&\leq &(1+\epsilon )\cdot {\bf B(}T)=(1+\epsilon )\cdot {\bf A}^{*\Delta
inj}(T)\text{ .}
\end{eqnarray*}
Since further 
\[
{\frak A}^{*}(Y^\infty ,X)\stackrel{1}{=}{\frak A}^{inj*}(Y^\infty ,X)\text{%
, } 
\]
(cf. \cite{df}, 20.12.), conjugation leads to the remaining incluson ${\frak %
A}^{inj*\Delta }\stackrel{1}{\subseteq }{\frak A}^{*\Delta inj}$, which
implies the equality 
\[
{\bf A}^{inj*\Delta }(T)={\bf A}^{inj\min }(T)={\bf A}^{*\Delta inj}(T) 
\]
for all $T\in {\frak F}$, so that in particular ${\frak A}^{inj*\Delta }%
\stackrel{\frak F}{=}{\frak A}^{inj\min }\stackrel{\frak F}{=}{\frak B}$ is
totally accessible. Since in general the inclusions ${\frak A}^{inj\min }%
\stackrel{1}{\subseteq }{\frak A}^{\min inj}$ (cf. \cite{df}, 25.11.) and $%
{\frak A}^{\min }\stackrel{1}{\subseteq }{\frak A}^{*\Delta }$ (cf. \cite
{oe3}) are satisfied, the proof is finished. $\blacksquare $

Due to the existence of Banach spaces without the metric approximation
property, we have ${\frak L}^{*\Delta inj}\stackrel{1}{=}{\frak I}^{\Delta
inj}\stackrel{1}{=}{\frak L}\stackrel{1}{\neq }{\frak I}^\Delta \stackrel{1}{%
=}{\frak L}^{*\Delta }$ and ${\frak K}\stackrel{1}{=}\overline{{\frak F}}%
^{inj}\stackrel{1}{=}{\frak L}^{\min inj}\stackrel{1}{\neq }{\frak L}%
^{inj\min }\stackrel{1}{=}\overline{{\frak F}}$, so that in general we
cannot transfer the previous Proposition to operators with infinite
dimensional range. What about quasi--Banach ideals which are not normed? As
the proof shows, the assumption $p=1$ is essential. But we even can say
more: The statement is false{\it \ }if we only assume the case $0<p<1$! To
see this, consider the {\it injective} $\frac 12$ --Banach ideal ${\frak A}:=%
{\frak P}_2\circ {\frak P}_2\stackrel{1}{=}{\frak L}_2\circ {\frak N}$ (cf. 
\cite{jo}, \cite{oe4}, \cite{p2}). Being a trace ideal, ${\frak A}$ cannot
be normed. Since the self--adjoint Banach ideal ${\frak P}_2$ is accessible,
the quotient formula implies that ${\frak A}^{inj*}\stackrel{1}{=}{\frak A}%
^{*}\stackrel{1}{=}{\frak P}_2^{*}\circ {\frak P}_2^{-1}\stackrel{1}{=}%
{\frak P}_2\circ {\frak P}_2^{-1}\stackrel{1}{=}{\frak L}$ (cf. \cite{df},
25.7.). In particular the ${\frak A}^{*}-LRP$ is valid, and we obtain $%
{\frak A}^{inj*\Delta }\stackrel{1}{=}{\frak A}^{*\Delta }\stackrel{1}{=}%
{\frak L}^\Delta \stackrel{1}{=}{\frak I}$. On the other hand ${\frak A}%
^{*\Delta inj}\stackrel{1}{=}{\frak I}^{inj}\stackrel{1}{=}{\frak P}_1$. The
assumption ${\frak I}\stackrel{\frak F}{=}{\frak P}_1$ would imply the
(global) equality ${\frak I}^\Delta \stackrel{1}{=}{\frak P}_1^\Delta 
\stackrel{1}{=}{\frak L}_\infty $ which is a contradiction\footnote{%
Note, that these considerations even show that ${\frak A}^{inj**}\stackrel%
{\frak F}{\neq }{\frak A}^{**inj}$.}, since ${\frak L}_2\stackrel{1}{=}%
{\frak D}_2^\Delta \stackrel{1}{\subseteq }{\frak I}^\Delta $.

However, there exists an additional sufficient condition which allows the
extension of the previous Proposition to operators between infinite
dimensional Banach spaces, namely the property (I) of the product ideal $%
{\frak A}^{*}\circ {\frak L}_\infty $:

\begin{theorem}
Let $({\frak A},{\bf A})$ be a maximal Banach ideal so that the ${\frak A}%
^{*}-LRP$ is satisfied. Then 
\begin{equation}
{\frak A}^{*\Delta inj}\stackrel{1}{\subseteq }({\frak A}^{*\Delta inj})^{dd}
\tag{$*$}
\end{equation}
If in addition, ${\frak A}^{*}\circ {\frak L}_\infty $ has the property (I),
then 
\begin{equation}
({\frak A}^{*\Delta inj})^{dd}\stackrel{1}{=}{\frak A}^{*\Delta inj}%
\stackrel{1}{=}{\frak P}_1\circ ({\frak A}^{*})^{-1}\stackrel{1}{=}{\frak A}%
^{inj*\Delta }\stackrel{1}{=}({\frak A}^{inj*\Delta })^{dd}  \tag{$**$}
\end{equation}
\end{theorem}

{\sc Proof:} First, let the ${\frak A}^{*}-LRP$ be satisfied. Let $T\in 
{\frak A}^{*\Delta inj}(X,Y)$ be given and $X$, $Y$ be arbitrary Banach
spaces. Due to Corollary 3.1 and the assumed validity of the ${\frak A}%
^{*}-LRP$, it follows that $J_Y^{\prime \prime }T^{\prime \prime
}=(J_YT)^{\prime \prime }\in $ ${\frak A}^{*\Delta }(X^{\prime \prime
},(Y^\infty )^{\prime \prime })$ and 
\[
{\bf A}^{*\Delta }(J_Y^{\prime \prime }T^{\prime \prime })\leq {\bf A}%
^{*\Delta }(J_YT)=({\bf A}^{*\Delta })^{inj}(T)\text{.} 
\]
Since $J_Y^{\prime \prime }:Y^{\prime \prime }\stackrel{1}{\hookrightarrow }%
(Y^\infty )^{\prime \prime }$ is an isometric embedding (cf. \cite{p1},
B.3.9.), the metric extension property of $(Y^{\prime \prime })^\infty $
implies the existence of an operator $\widetilde{J}\in {\frak L}((Y^\infty
)^{\prime \prime },(Y^{\prime \prime })^\infty )$ so that $J_{Y^{\prime
\prime }}=\widetilde{J}J_Y^{\prime \prime }$ and $\Vert \widetilde{J}\Vert
=1 $. Hence, $T^{\prime \prime }\in ({\frak A}^{*\Delta })^{inj}(X^{\prime
\prime },Y^{\prime \prime })$ and 
\[
({\bf A}^{*\Delta })^{inj}(T^{\prime \prime })\leq ({\bf A}^{*\Delta
})^{inj}(T)\text{,} 
\]
which implies the inclusion $(*)$. To prove $(**)$, note, that the second
isometric identity already has been proven in this paper (see Proposition
2.2). Recalling that always 
\[
{\frak A}^{inj*\Delta }\stackrel{1}{\subseteq }({\frak A}^{*\Delta })^{inj}%
\text{,} 
\]
we only have to prove the inclusion 
\[
({\frak A}^{*\Delta inj})^{dd}\stackrel{1}{\subseteq }{\frak A}^{inj*\Delta
} 
\]
-- given the property (I) of ${\frak A}^{*}\circ {\frak L}_\infty $. To this
end, let $T$ $\in ({\frak A}^{*\Delta inj})^{dd}(X,Y)$ be given, with
arbitrarily chosen Banach spaces $X$ and $Y$, and put $({\frak B},{\bf B}):=(%
{\frak A}^{inj},{\bf A}^{inj})$. Since ${\frak B}^{*\Delta }$ is regular
(see Proposition 2.1), we only have to show that $j_YT\in {\frak B}^{*\Delta
}(X,Y^{\prime \prime })$ and 
\[
{\bf B}^{*\Delta }(j_YT)\leq ({\bf A}^{*\Delta })^{inj}(T^{\prime \prime })%
\text{.} 
\]
So, let $L\in {\frak F}(Y^{\prime \prime },X)$ be an arbitrary finite rank
operator -- considered as an element of ${\frak B}^{*}(Y^{\prime \prime },X)$%
. Due to the assumed property (I) of ${\frak A}^{*}\circ {\frak L}_\infty $,
Theorem 3.3 shows us the existence of a {\it finite rank} operator $V\in 
{\frak F}((Y^{\prime \prime })^\infty ,X^{\prime \prime })$ so that $%
j_XL=VJ_{Y^{\prime \prime }}$ and 
\[
{\bf A}^{*}(V)\leq (1+\epsilon )\cdot {\bf B}^{*}(L)\text{.} 
\]
Hence, 
\begin{eqnarray*}
&\mid &tr(j_YTL)\mid =\mid tr(T^{\prime \prime }j_XL)\mid =\mid tr(T^{\prime
\prime }VJ_{Y^{\prime \prime }})\mid =\mid tr(J_{Y^{\prime \prime
}}T^{\prime \prime }V)\mid \\
&\leq &{\bf A}^{*\Delta }(J_{Y^{\prime \prime }}T^{\prime \prime })\cdot 
{\bf A}^{*}(V) \\
&\leq &(1+\epsilon )\cdot ({\bf A}^{*\Delta })^{inj}(T^{\prime \prime
})\cdot {\bf B}^{*}(L)\text{,}
\end{eqnarray*}
which implies that $j_YT\in {\frak B}^{*\Delta }(X,Y^{\prime \prime })$ and $%
{\bf B}^{*\Delta }(j_YT)\leq ({\bf A}^{*\Delta })^{inj}(T^{\prime \prime })$%
. Summing up all the previous steps in our proof, we have shown that 
\[
{\frak A}^{*\Delta inj}\stackrel{1}{=}{\frak A}^{inj*\Delta }\stackrel{1}{=}(%
{\frak A}^{*\Delta inj})^{dd} 
\]
which obviously implies $(**)$, and the proof is finished.$\blacksquare $

\begin{corollary}
Let $({\frak A},{\bf A})$ be a maximal and left--accessible Banach ideal so
that ${\frak A}^{*}\circ {\frak L}_\infty $ has the property (I). Then both, 
${\frak A}^{inj}$ and $({\frak A}^{inj})^{*}$ are totally accessible.
\end{corollary}

{\sc Proof:} Since ${\frak A}$ is left--accessible, Proposition 2.2 and the
previous statement imply that 
\[
{\frak A}^{inj}\stackrel{1}{\subseteq }{\frak P}_1\circ ({\frak A}^{*})^{-1}%
\stackrel{1}{=}({\frak A}^{*\Delta })^{inj}\stackrel{1}{=}({\frak A}%
^{inj})^{*\Delta }\stackrel{1}{\subseteq }{\frak A}^{inj}\text{,}
\]
and it follows that 
\[
{\frak B}^{*}\stackrel{1}{=}{\frak B}^\Delta 
\]
where we have put ${\frak B}:=({\frak A}^{inj})^{*}$. Hence, 
\cite[Theorem 3.1.]{oe4} finishes the proof.$\blacksquare $

As the careful reader might guess, Theorem 3.4 and Corollary 3.4 imply a lot
of interesting consequences. Let us note the most important one:

\begin{theorem}
Let $({\frak A},{\bf A})$ be a maximal Banach ideal so that ${\frak A}%
^{*}\circ {\frak L}_\infty $ has the property (I). If space$({\frak A})$
contains a Banach space $X_0$ so that $X_0$ has the bounded approximation
property but $X_0^{\prime \prime }$ has not, then the ${\frak A}^{*}-LRP$
cannot be satisfied.
\end{theorem}

{\sc Proof:} Assume, that the statement is false and hence the ${\frak A}%
^{*}-LRP$ is satisfied. Since $X_0$ has the bounded approximation property, $%
Id_{X_0}\in {\frak I}^\Delta (X_0,X_0)$ and $c:={\bf I}^\Delta
(Id_{X_0})<\infty $. By definition of ${\frak I}^\Delta $ and of the adjoint 
${\frak A}^{*}$, one immediately derives the inclusion 
\[
{\frak A}\circ {\frak I}^\Delta \circ {\frak A}^{*}\stackrel{1}{\subseteq }%
{\frak I}\text{,} 
\]
so that in particular 
\[
{\frak A}(X_0,X_0)\subseteq {\frak I}\circ ({\frak A}^{*})^{-1}(X_0,X_0)%
\stackrel{1}{\subseteq }{\frak P}_1\circ ({\frak A}^{*})^{-1}(X_0,X_0)%
\stackrel{1}{=}{\frak A}^{*\Delta inj}(X_0,X_0) 
\]
and 
\[
{\bf A}^{*\Delta inj}(Id_{X_0})\leq c\cdot {\bf A}(Id_{X_0})\text{.} 
\]
Hence, due to the assumed property (I) of ${\frak A}^{*}\circ {\frak L}%
_\infty $, Theorem 3.4 implies that even $X_0^{\prime \prime }\in $ space$(%
{\frak A}^{inj*\Delta })$ and 
\[
{\bf A}^{inj*\Delta }(Id_{X_0^{\prime \prime }})\stackrel{1}{=}{\bf A}%
^{inj*\Delta }(Id_{X_0}^{\prime \prime })\leq c\cdot {\bf A}(Id_{X_0})\text{.%
} 
\]
But this would imply that $X_0^{\prime \prime }\in $ space$({\frak I}^\Delta
)$, leading to the conclusion that $X_0^{\prime \prime }$ would have the
bounded approximation property -- with constant $c\cdot {\bf A}(Id_{X_0})$,
which is a contradiction.$\blacksquare $

Now, the reader may ask for explicite examples for such maximal Banach
ideals. To this end, note again that ${\frak A}^{*}\circ {\frak L}_\infty $
has the property (I), if ${\frak A}^{*}$ contains ${\frak L}_2$ as a factor.
Since ${\frak A}^{*}$ is a Banach ideal, we therefore have to look for
maximal {\it Banach} ideals of type ${\frak B}\circ {\frak L}_2\circ {\frak C%
}$. A first investigation of geometrical properties of such product ideals
was given in \cite{oe4}. Unfortunately, we cannot present explicite
sufficient criteria which show the existence of (an equivalent) ideal norm
on product ideals. It seems to be much more easier to show that a certain
product ideal cannot be a normed one by using arguments which involve trace
ideals and the ideal of nuclear operators (the smallest Banach ideal).
However, let us turn to the following section.

\section{Applications}

Among other things, we will see in this section how deep the properties (I)
and (S) reflect the local structure of operator ideals. A first example
considers the question of Defant and Floret (see \cite{df}, 21.12) whether $%
{\frak L}_\infty $ is totally accessible or not. We are able to show that $%
{\frak L}_\infty $ {\it is not totally accessible}, and the idea of the
proof is the following: Assuming the opposite, leads to the property (I) for
a suitable class of quasi--Banach ideals of type ${\frak A}^{*}\circ {\frak L%
}_\infty $. On the other hand, there exists a well known left--accessible
candidate ${\frak A}$ so that $({\frak A}^{inj})^{*}$ is not totally
accessible -- a contradiction to Corollary 3.4. To prepare the steps
carefully, we first state a fact which is of its own interest:

\begin{lemma}
Let $({\frak A},{\bf A})$ and $({\frak B},{\bf B})$ be arbitrary
quasi--Banach ideals so that

\begin{enumerate}
\item[(i)]  ${\frak A}\circ {\frak B}$ has the property (S)

\item[(ii)]  ${\frak B}$ is totally accessible.
\end{enumerate}

Then ${\frak A}\circ {\frak B}$ is left--accessible and has the property (I).
\end{lemma}

{\sc Proof:} Let $X$, $Y$ be arbitrary Banach spaces and $L\in {\frak F}%
(X,Y) $ an arbitrary finite rank operator. Given $\epsilon >0$, there exists
a Banach space $Z$ and operators $A\in {\frak A}(Z,Y)$, $B\in {\frak B}(X,Z)$
so that $L=AB$ and 
\[
{\bf A}(A)\cdot {\bf B}(B)\leq (1+\epsilon )\cdot {\bf A\circ B}(L)\text{.} 
\]
Due to the property (S) of ${\frak A}\circ {\frak B}$, we may assume that $B$
is of finite rank. Hence, since ${\frak B}$ is totally accessible, there
exist $K\in $ COFIN$(X)$, $E\in $ FIN$(Z)$ and an operator $\Gamma \in 
{\frak L}(X\diagup K,E)$ so that $B=J_E^Z\Gamma Q_K^X$ and 
\[
{\bf B}(\Gamma )\leq (1+\epsilon )\cdot {\bf B}(B)\text{.} 
\]
Therefore, $L=A_0\Gamma Q_K^X$ where $A_0:=AJ_E^Z\in {\frak F}(E,Y)$ and 
\[
{\bf A\circ B(}A_0\Gamma )\leq {\bf A}(A_0)\cdot {\bf B}(\Gamma )\leq
(1+\epsilon )^2\cdot {\bf A\circ B}(L)\text{,} 
\]
and the claim follows.$\blacksquare $

Obviously, similar arguments allow a transfer of property (S) to property
(I), and we obtain the ''(I)--version'':

\begin{lemma}
Let $({\frak A},{\bf A})$ and $({\frak B},{\bf B})$ be arbitrary
quasi--Banach ideals so that

\begin{enumerate}
\item[(i)]  ${\frak A}\circ {\frak B}$ has the property (I)

\item[(ii)]  ${\frak A}$ is totally accessible.
\end{enumerate}

Then ${\frak A}\circ {\frak B}$ has the property (S) and is
right--accessible.
\end{lemma}

Associating these results with Proposition 2.3, we immediately obtain (with
the help of a factor diagram) a quite useful result\footnote{%
Notice, that Proposition 21.4. in \cite{df} is a Corollary of this result.}:

\begin{proposition}
Let $({\frak A},{\bf A})$ and $({\frak B},{\bf B})$ be arbitrary
quasi--Banach ideals so that one of the following properties hold:

\begin{enumerate}
\item[(1)]  ${\frak A}\circ {\frak B}$ has the property (I), ${\frak A}$ is
totally accessible and ${\frak B}$ is left--acessible

\item[(2)]  ${\frak A}\circ {\frak B}$ has the property (S), ${\frak A}$ is
right--accessible and ${\frak B}$ is totally accessible,
\end{enumerate}

then ${\frak A}\circ {\frak B}$ has the property (I) as well as the property
(S), and $({\frak A}\circ {\frak B})^{reg}$ is totally accessible.
\end{proposition}

Now, we are well prepared to investigate the total accessibility of ${\frak L%
}_\infty $:

\begin{theorem}
The maximal Banach ideals ${\frak L}_\infty \thicksim g_\infty $ and ${\frak %
L}_1\thicksim w_1$ are not totally accessible.
\end{theorem}

{\sc Proof:} Since ${\frak L}_1\stackrel{1}{=}{\frak L}_\infty ^d$, we only
have to prove the claim for ${\frak L}_\infty $. Assume the opposite.
Consider the maximal Banach ideal ${\frak A}:={\frak L}_1\stackrel{1}{=}(%
{\frak P}_1^d)^{*}$. Since 
\[
({\frak A}^{*})^{sur}\stackrel{1}{=}({\frak P}_1^d)^{sur}\stackrel{1}{=}(%
{\frak P}_1^{inj})^d\stackrel{1}{=}{\frak P}_1^d\stackrel{1}{=}{\frak A}^{*}%
\text{,} 
\]
it follows that ${\frak A}^{*}$ is surjective, so that ${\frak A}^{*}\circ 
{\frak L}_\infty $ has the property (S). Due to Lemma 4.1, the assumed total
accessibility of ${\frak L}_\infty $ even leads to the property (I) of $%
{\frak A}^{*}\circ {\frak L}_\infty $, and Corollary 3.4 implies that $(%
{\frak A}^{inj})^{*}\stackrel{1}{=}({\frak L}_1^{inj})^{*}$ is totally
accessible. On the other hand, \cite[Corollary 21.6.2]{df} tells us that the
adjoint of ${\frak L}_1^{inj}$ cannot be totally accessible (because of the
existence of subspaces of $l_1$ without the approximation property), and we
obtain a contradiction.$\blacksquare $

\begin{corollary}
${\frak L}_\infty \circ {\frak L}_\infty $ neither has property (I) nor
property (S) and is not regular. In particular, ${\frak L}_\infty \stackrel{1%
}{\neq }{\frak L}_\infty \circ {\frak L}_\infty $.
\end{corollary}

{\sc Proof:} First, assume that ${\frak L}_\infty \circ {\frak L}%
_\infty $ has property (I). Then, Proposition 3.1 implies that 
\[
{\frak L}_\infty ^{-1}\circ {\frak L}_\infty ^\Delta \stackrel{1}{=}({\frak L%
}_\infty \circ {\frak L}_\infty )^\Delta \text{.}
\]
Since ${\frak P}_1$ is right--accessible, it follows that 
\[
{\frak L}_\infty \circ {\frak P}_1\stackrel{1}{=}{\frak P}_1^{*}\circ {\frak %
P}_1\stackrel{1}{\subseteq }{\frak I}\stackrel{1}{=}{\frak L}^\Delta 
\stackrel{1}{\subseteq }{\frak L}_\infty ^\Delta \text{,}
\]
and hence ${\frak P}_1\stackrel{1}{\subseteq }({\frak L}_\infty \circ {\frak %
L}_\infty )^\Delta $. But this inclusion further implies that 
\[
{\frak P}_1^d\stackrel{1}{\subseteq }({\frak L}_\infty \circ {\frak L}%
_\infty )^{d\Delta }\stackrel{1}{=}(({\frak L}_\infty \circ {\frak L}_\infty
)^{reg})^{d\Delta }\stackrel{1}{=}{\frak L}_\infty ^{d\Delta }\stackrel{1}{=}%
{\frak L}_1^\Delta \text{,}
\]
(since ${\frak L}_\infty \stackrel{1}{=}({\frak L}_\infty \circ {\frak L}%
_\infty )^{reg}$) and we obtain the contradiction ${\frak L}_1^{*}\stackrel{1%
}{=}{\frak L}_1^\Delta $ (since ${\frak L}_1$ is not totally accessible).
Using similar arguments, the assumption of property (S) of ${\frak L}_\infty
\circ {\frak L}_\infty $ also leads to a contradiction. Obviously, ${\frak L}%
\circ {\frak L}_\infty \stackrel{1}{=}{\frak L}_\infty $ has the property
(S), and the proof is finished.$\blacksquare $

Given two ultrastable quasi--Banach ideals $({\frak A},{\bf A})$ and $(%
{\frak B},{\bf B})$, we know that the product ideal ${\frak B}\circ {\frak A}
$ is also ultrastable (see \cite{de}). If ${\frak B}\circ {\frak A}$ is
normed, Theorem 2.1 further implies that 
\[
({\frak B}\circ {\frak A})^{**}\stackrel{1}{=}({\frak B}\circ {\frak A}%
)^{\max }\stackrel{1}{=}({\frak B}\circ {\frak A})^{reg}\text{,} 
\]
and, if {\it in addition} $({\frak A},{\bf A})$ is {\it injective}, then 
\cite[Corollary 2.1.]{oe3} even implies that 
\[
({\frak B}\circ {\frak A})^{**}\stackrel{1}{=}({\frak B}\circ {\frak A}%
)^{reg}\stackrel{1}{=}{\frak B}^{reg}\circ {\frak A}\text{ .} 
\]
Bearing this situation in mind, the ${\frak B}-LRP$ leads to a further
surprising result which has rich consequences:

\begin{lemma}
Let $({\frak A},{\bf A})$ and $({\frak B},{\bf B})$ be quasi--Banach ideals,
so that ${\frak L}_2$ is a left--factor of ${\frak A}$ and the ${\frak B}-LRP
$ is satisfied. Then, 
\[
{\frak B}\circ {\frak A}\stackrel{\frak F}{=}({\frak B}\circ {\frak A})^{reg}%
\stackrel{1}{=}{\frak B}^{reg}\circ {\frak A}\text{ .}
\]
\end{lemma}

{\sc Proof:} Let $X$, $Y$ be arbitrary Banach spaces, $\epsilon >0$ and $%
L\in {\frak F}(X,Y)$ be an arbitrary {\it finite rank} operator. Let $(%
{\frak C},{\bf C})$ be a quasi--Banach ideal so that ${\frak A}\stackrel{1}{=%
}{\frak L}_2\circ {\frak C}$. Due to the property (S) of the product ideal $%
{\frak B}\circ {\frak A}$ and the property (I) of (the injective product) $%
{\frak L}_2\circ {\frak C}$, we may write the finite rank operator $j_YL$ as
the composition 
\[
j_YL=B\Lambda C\text{,} 
\]
where $C\in {\frak C}(X,U)$, $\Lambda \in {\frak F}(U,V)$, $B\in {\frak B}%
(V,Y^{\prime \prime })$ and $U$, $V$ are Banach spaces so that 
\[
{\bf B}(B){\bf \cdot L}_2(\Lambda )\cdot {\bf C}(C)<(1+\epsilon )^2\cdot (%
{\bf B\circ A)}^{reg}(L)\text{.} 
\]
Since ${\frak L}_2$ is totally accessible, there exists a finite dimensional
subspace $F$ of $V$ and an operator $\Lambda _0\in {\frak L}(U,F)$ so that $%
\Lambda =J_F^V\Lambda _0$ and ${\bf L}_2(\Lambda _0)<(1+\epsilon )\cdot {\bf %
L}_2(\Lambda )$. Hence, we now may apply the assumed ${\frak B}-LRP$ to the
operator $BJ_F^V\in {\frak L}(F,Y^{\prime \prime })$, and it therefore
follows the existence of an operator $B_0\in {\frak L}(F,Y)$ so that 
\[
{\bf B}(B_0)\leq (1+\epsilon )\cdot {\bf B}(B) 
\]
and 
\[
BJ_F^Vv=j_YB_0v\text{\ \ for all }v\in (BJ_F^V)^{-1}(j_Y(Y))\text{.} 
\]
Let $x\in X$ be given. Then $v=$ $\Lambda _0Cx\in (BJ_F^V)^{-1}(j_Y(Y))$,
which implies that 
\[
j_YLx=B\Lambda Cx=BJ_F^Yv=j_YB_0v=j_YB_0\Lambda _0Cx\text{.} 
\]
Since $x\in X$ was chosen arbitrary, we therefore obtain that $L=B_0\Lambda
_0C$ and 
\begin{eqnarray*}
{\bf B}\circ {\bf A}(L) &=&{\bf B}\circ {\bf L}_2\circ {\bf C}(L)\leq {\bf B}%
(B_0)\cdot {\bf L}_2(\Lambda _0)\cdot {\bf C}(C) \\
&\leq &(1+\epsilon )^2\cdot {\bf B}(B)\cdot {\bf L}_2(\Lambda )\cdot {\bf C}%
(C) \\
&\leq &(1+\epsilon )^4\cdot ({\bf B\circ A)}^{reg}(L)\text{,}
\end{eqnarray*}
and the proof is finished.$\blacksquare $

For completion, let us note the following fact:

\begin{proposition}
Let $({\frak A},{\bf A})$ be an arbitrary quasi--Banach ideal. Then, ${\frak %
A}^{*}\circ {\frak L}_2\circ {\frak A}$ cannot be a $1$-- Banach ideal.
\end{proposition}

{\sc Proof:} Put ${\frak B}:={\frak L}_2\circ {\frak A}$. Then both, ${\frak %
B}$ and ${\frak B}\circ \overline{{\frak F}}$ are injective (see \cite{oe4}%
). Assume, that the statement is false. The previous considerations imply
that the regular hull $({\frak A}^{*}\circ {\frak B}\circ \overline{{\frak F}%
})^{reg}$ coincides isometrically with ${\frak A}^{*}\circ {\frak B}\circ 
\overline{{\frak F}}$. Using exactly the same technique as presented in the
proof of Proposition 8.3. in \cite{b} (a factorization through Banach spaces
of type $l_\infty ^2(Z_1,Z_2)$ consisting of elements $(z_1,z_2)\in
Z_1\times Z_2$ so that $\Vert (z_1,z_2)\Vert _\infty :=\max (\Vert z_1\Vert
,\Vert z_2\Vert )<\infty $) shows, that the assumed existence of an ideal
norm on ${\frak A}^{*}\circ {\frak B}$ even implies the existence of an
ideal norm on the smaller ideal ${\frak A}^{*}\circ {\frak B}\circ \overline{%
{\frak F}}$. Since ${\frak N}$, the ideal of the nuclear operators is the
smallest Banach ideal, the right--accessibility of ${\frak B}$ therefore
implies that ${\frak N}\stackrel{1}{\subseteq }{\frak A}^{*}\circ {\frak B}%
\circ \overline{{\frak F}}\stackrel{1}{\subseteq }{\frak B}^{*}\circ {\frak B%
}\circ \overline{{\frak F}}\stackrel{1}{\subseteq }{\frak N}$, and we obtain
that 
\[
{\frak N}\stackrel{1}{=}{\frak A}^{*}\circ {\frak B}\circ \overline{{\frak F}%
}\stackrel{1}{=}({\frak A}^{*}\circ {\frak B}\circ \overline{{\frak F}}%
)^{reg}\stackrel{1}{=}{\frak N}^{reg}\stackrel{1}{=}{\frak N}^d 
\]
-- a contradiction (cf. \cite[Proposition 16.8.]{df}).$\blacksquare $

Next, we again turn our attention to candidates $({\frak A},{\bf A})$ which
do not satisfy the ${\frak A}-LRP$. Although it seems, that property (I) and
property (S) do not play the fundamental part in the next statement, the
proof opens the reader`s eyes:

\begin{theorem}
Let $({\frak A},{\bf A)}$ be a maximal Banach ideal so that the injective
hull of $({\frak A}\circ {\frak L}_2)^{\max }$ is not totally accessible,
then the $({\frak A}^{inj})^{*}-LRP$ (and hence the ${\frak A}^{*}-LRP$)
cannot be satisfied.
\end{theorem}

{\sc Proof:} Assume that the $({\frak A}^{inj})^{*}-LRP$ holds. Then $(%
{\frak A}^{inj})^{*\Delta }$ is totally accessible -- due to Proposition
3.2. Since ${\frak L}_2$ is injective, ${\frak A}_0:=({\frak A}%
^{inj})^{*\Delta }\circ {\frak L}_2$ has the property (I), and the total
accessibility of $({\frak A}^{inj})^{*\Delta }$ together with the
left--accessibility of ${\frak L}_2$ further implies that ${\frak A}_0$ is
totally accessible (Proposition 4.1). Since ${\frak L}_2$ is a factor of $%
{\frak A}_0$, ${\frak L}_\infty \circ {\frak A}_0$ has the property (S), and
Proposition 4.1 even implies that $({\frak L}_\infty \circ {\frak A}%
_0)^{reg} $ is totally accessible. Because of the metric approximation
property of Hilbert spaces and spaces of type $L_\infty $, $({\frak L}%
_\infty \circ {\frak A}_0)^{reg}\stackrel{1}{=}({\frak L}_\infty \circ 
{\frak A}^{inj}\circ {\frak L}_2)^{reg}$, so that in particular ${\frak L}%
_\infty \circ {\frak A}^{inj}\circ {\frak L}_2$ is totally accessible.
Therefore, we obtain the total accessibility of the injective hull of $%
{\frak L}_\infty \circ {\frak A}^{inj}\circ {\frak L}_2$ which (by
definition) equals isometrically $({\frak A}^{inj}\circ {\frak L}_2)^{inj}%
\stackrel{1}{=}({\frak A}\circ {\frak L}_2)^{inj}$. Hence, $(({\frak A}\circ 
{\frak L}_2)^{\max })^{inj}\stackrel{1}{=}(({\frak A}\circ {\frak L}%
_2)^{reg})^{inj}\stackrel{1}{=}({\frak A}\circ {\frak L}_2)^{inj}$ must be
totally accessible.$\blacksquare $

\begin{corollary}
Let $({\frak B},{\bf B)}$ be a maximal Banach ideal so that ${\frak A}:=%
{\frak B}\circ {\frak L}_2$ is normed. Then, ${\frak A}$ is a maximal Banach
ideal, and the following statements are equivalent:

\begin{enumerate}
\item[(i)]  The $({\frak A}^{inj})^{*}-LRP$ is satisfied

\item[(ii)]  ${\frak A}^{inj}$ is totally accessible.
\end{enumerate}
\end{corollary}

{\sc Proof:} Only the inclusion (i)$\Longrightarrow $(ii) is not trivial.
Since 
\[
{\frak A}^{**}\stackrel{1}{=}{\frak A}^{\max }\stackrel{1}{=}({\frak B}\circ 
{\frak L}_2)^{reg}\stackrel{1}{=}{\frak B}\circ {\frak L}_2\stackrel{1}{=}%
{\frak A}\text{,} 
\]
${\frak A}$ is a maximal Banach ideal, and we therefore may apply Theorem
4.2 to ${\frak A}$. Since ${\frak L}_2\stackrel{1}{=}{\frak L}_2\circ {\frak %
L}_2$, assumption (i) therefore implies that ${\frak A}^{inj}\stackrel{1}{=}(%
{\frak A}\circ {\frak L}_2)^{inj}\stackrel{1}{=}(({\frak A}\circ {\frak L}%
_2)^{\max })^{inj}$ is totally accessible, and the Corollary is proven.$%
\blacksquare $

The careful reader now may (and should) ask whether there exist operator
ideals ${\frak A}$ so that ${\frak A}\circ {\frak L}_2$ is not injective.
Indeed, we will show that this is the case -- in contrast to ideals of type $%
{\frak L}_2\circ {\frak A}$ which always are injective{\em \ }(see \cite{oe4}%
, Lemma 5.1). To this end, we need the help of some ''exotic'' Banach
spaces: G.T. spaces. Recall that a Banach space $X$ is called a G.T. space
(a space which satisfies Grothendieck`s Theorem) if 
\[
{\frak L}(X,l_2)={\frak P}_1(X,l_2) 
\]
Details and further informations about these Banach spaces are listed in 
\cite{df} and \cite{pi}. We now will work with the famous Pisier space $P$
which is a G.T. space without the approximation property (see \cite{pi},
Theorem 10.6.). By the ${\frak L}_p$ -- Local Technique Lemma for Operator
Ideals (see \cite{df}, 23.1.), it follows that ${\frak L}_2(P,\cdot
)\subseteq {\frak P}_1(P,\cdot )$ which is equivalent to $Id_P\in {\frak L}%
_2^{-1}\circ {\frak P}_1\stackrel{1}{=}({\frak L}_\infty \circ {\frak L}%
_2)^{*}$ (since ${\frak L}_2$ is totally accessible). Assuming now that $%
{\frak L}_\infty \circ {\frak L}_2$ is injective, would imply the
contradiction $Id_P\in ({\frak L}_\infty \circ {\frak L}_2)^{*}\stackrel{1}{=%
}(({\frak L}_\infty \circ {\frak L}_2)^{inj})^{*}\stackrel{1}{=}{\frak L}%
_2^{*}\stackrel{1}{=}{\frak L}_2^\Delta \stackrel{1}{\subseteq }{\frak I}%
^\Delta $. Hence, ${\frak L}_\infty \circ {\frak L}_2$ is not injective. But
even more holds:

\begin{proposition}
Let $({\frak A},{\bf A)}$ be a maximal Banach ideal which contains ${\frak L}%
_2$ as a factor. Then ${\frak L}_\infty \circ {\frak A}$ is not injective.
\end{proposition}

{\sc Proof:} Put ${\frak B}_0:=({\frak L}_\infty \circ {\frak L}_2)^{*}$.
Then ${\frak B}_0\subseteq ({\frak L}_\infty \circ {\frak A})^{*}$. Assuming
the injectivity of ${\frak L}_\infty \circ {\frak A}$ leads to the inclusion 
${\frak B}_0\subseteq ({\frak A}^{inj})^{*}$, and we obtain 
\[
{\frak B}_0\circ {\frak A}\subseteq ({\frak A}^{inj})^{*}\circ {\frak A}%
\stackrel{1}{\subseteq }({\frak A}^{inj})^{*}\circ {\frak A}^{inj}\stackrel{1%
}{\subseteq }{\frak I} 
\]
(since ${\frak A}^{inj}$ always is right--accessible), which implies that $%
Id_P\in {\frak B}_0\subseteq {\frak I}\circ {\frak A}^{-1}\stackrel{1}{%
\subseteq }{\frak P}_1\circ {\frak A}^{-1}\stackrel{1}{=}({\frak A}^\Delta
)^{inj}$. Since ${\frak L}_2$ is a factor of ${\frak A}$, ${\frak A}\circ 
{\frak L}_\infty $ has the property (I), and the proof of Theorem 3.4
implies that $Id_p\in {\frak B}_0\stackrel{1}{=}{\frak B}_0^{dd}\subseteq ((%
{\frak A}^\Delta )^{inj})^{dd}\stackrel{1}{\subseteq }({\frak A}%
^{*inj})^{*\Delta }\stackrel{1}{\subseteq }{\frak I}^\Delta $ which is a
contradiction.$\blacksquare $

To round off these interesting considerations, we next prove a quotient
version of Grothendieck`s Theorem:

\begin{proposition}
Let ${\frak B}_0:=({\frak L}_\infty \circ {\frak L}_2)^{*}$. Then ${\frak L}%
_1\varsubsetneq {\frak B}_0$, and space$({\frak B}_0)$ contains Banach
spaces without the approximation property. Moreover, 
\[
{\frak L}_2\stackrel{1}{=}{\frak P}_1\circ {\frak B}_0^{-1}\text{.}
\]
\end{proposition}

{\sc Proof:} The inclusion $\stackrel{1}{\subseteq }$ already has been
shown. To see the other inclusion, note that ${\frak D}_2\stackrel{1}{=}%
{\frak L}_2^{*}\stackrel{1}{\subseteq }{\frak B}_0$. Since ${\frak L}_2$ is
injective, it therefore follows that ${\frak P}_1\circ {\frak B}_0^{-1}%
\stackrel{1}{=}({\frak B}_0^\Delta )^{inj}\stackrel{1}{\subseteq }({\frak D}%
_2^\Delta )^{inj}\stackrel{1}{=}{\frak L}_2$.$\blacksquare $

The situation completely changes, if we permute the factors ${\frak L}%
_\infty $ and ${\frak L}_2$ in the product ideal ${\frak L}_\infty \circ 
{\frak L}_2$, since: 
\begin{equation}
({\frak L}_2\circ {\frak L}_\infty )^{*}\stackrel{1}{=}{\frak P}_1\circ 
{\frak L}_2^{-1}\stackrel{1}{=}{\frak D}_2^{inj}\stackrel{1}{\subseteq }%
{\frak P}_2\text{.}  \tag{$*$}
\end{equation}

If $({\frak B},{\bf B)}$ is a quasi--Banach ideal so that ${\frak B}$ $%
\subseteq {\frak D}_2$, then we already know that ${\frak B}\circ {\frak L}%
_2 $ is a trace ideal and therefore cannot admit an (equivalent) ideal--{\it %
norm} (see \cite{jo}, 3.7.). What can we say if we only assume the existence
of one (suitable) Banach space $X_0$ so that ${\frak B}(\cdot ,X_0)\subseteq 
{\frak D}_2(\cdot ,X_0)$ ? In this case, the existence of an ideal--norm on
the product ideal ${\frak B}\circ {\frak L}_2$ a priori cannot be excluded,
and we will see that the property (I) implies a surprising connection
between the principle of local reflexivity for operator ideals and the
existence of such an ideal--norm. To prepare the right instruments, we need
the following statement\footnote{%
Note, that we do not assume the regularity of $({\frak B},{\bf B})$ in Lemma
4.4!}

\begin{lemma}
Let $({\frak B},{\bf B})$ be an arbitrary ultrastable quasi--Banach ideal so
that ${\frak B}\circ {\frak L}_2$ is normed. If the $({\frak B}\circ {\frak L%
}_2)^{**}-LRP$ is satisfied, then $({\frak B}\circ {\frak L}_2)^{**}\circ 
{\frak L}_\infty $ has the property (I) as well the property (S).
\end{lemma}

{\sc Proof:} Put ${\frak A}:=({\frak B}\circ {\frak L}_2)^{*}$. Due to the
assumptions on ${\frak B}$ and the product ideal ${\frak B}\circ {\frak L}_2$%
, it follows that 
\[
{\frak A}^{*}\stackrel{1}{=}({\frak B}\circ {\frak L}_2)^{**}\stackrel{1}{=}(%
{\frak B}\circ {\frak L}_2)^{\max }\stackrel{1}{=}({\frak B}\circ {\frak L}%
_2)^{reg} 
\]
is a maximal Banach ideal which even implies that 
\[
{\frak A}^{*}\stackrel{1}{=}(({\frak B}\circ {\frak L}_2)^{reg})^{dd}%
\stackrel{1}{=}({\frak B}\circ {\frak L}_2)^{dd}\text{.} 
\]
Let $X$, $Y$ be arbitrary Banach spaces and $T\in {\frak F}(X,Y)$ an
arbitrary finite rank operator. Given $\epsilon >0$, the definition of $%
{\frak L}_\infty $ implies the existence of a Borel--Radon measure $\mu $, a
Banach space $Z$ and operators $S_1\in {\frak L}(X,L_\infty (\mu ))$, $%
S_2\in {\frak L}(L_\infty (\mu ),Z^{\prime \prime })$, $R\in {\frak A}%
^{*}(Z,Y)\stackrel{1}{=}({\frak B}\circ {\frak L}_2)^{dd}(Z,Y)$ so that $%
j_YT $ $=R^{\prime \prime }S_2S_1$ and 
\[
{\bf B}\circ {\bf L}_2(R^{\prime \prime }S_2)\cdot \Vert S_1\Vert \leq {\bf B%
}\circ {\bf L}_2(R^{\prime \prime })\cdot \Vert S_2\Vert \cdot \Vert
S_1\Vert <(1+\epsilon )^2\cdot {\bf A}^{*}\circ {\bf L}_\infty (T)\text{.} 
\]
Since ${\frak L}_2$ is a factor of the product ideal ${\frak B}\circ {\frak L%
}_2$, we may copy the proof of \cite[Lemma 2.4.]{jo}, which even allows us
to {\it substitute} the operator $R^{\prime \prime }S_2$ through a finite
rank operator $V\in {\frak F}(L_\infty (\mu ),Y^{\prime \prime })$ so that $%
j_YT$ $=VS_1$ and 
\[
{\bf B}\circ {\bf L}_2(V)\cdot \Vert S_1\Vert <(1+\epsilon )^3\cdot {\bf A}%
^{*}\circ {\bf L}_\infty (T)\text{.} 
\]
Due to the metric approximation property of $L_\infty (\mu )^{\prime }$, we
then obtain a finite dimensional subspace $F$ of $L_\infty (\mu )$ and
operators $B\in {\frak L}(L_\infty (\mu ),F)$ and $W\in {\frak L}%
(F,Y^{\prime \prime })$ so that $V=WB$, $\Vert B\Vert \leq 1+\epsilon $ and $%
{\bf B}\circ {\bf L}_2(W)\leq {\bf B}\circ {\bf L}_2(V)$. Now, we proceed as
in the proof of the second part of Theorem 3.3, and the assumed ${\frak A}%
^{*}-LRP$ even implies the existence of an operator $W_0\in {\frak L}(F,Y)$
so that $T=W_0(BS_1)$ and 
\begin{eqnarray*}
{\bf A}^{*}(W_0)\cdot {\bf L}_\infty (BS_1) &\leq &(1+\epsilon )\cdot {\bf A}%
^{*}(W_0)\cdot \Vert S_1\Vert \leq (1+\epsilon )^2\cdot {\bf A}^{*}(W)\cdot
\Vert S_1\Vert \\
&\leq &(1+\epsilon )^2\cdot {\bf B}\circ {\bf L}_2(W)\cdot \Vert S_1\Vert
\leq (1+\epsilon )^5\cdot {\bf A}^{*}\circ {\bf L}_\infty (T)\text{,}
\end{eqnarray*}
and we have obtained the properties (I) and (S) of ${\frak A}^{*}\circ 
{\frak L}_\infty $.$\blacksquare $

\begin{theorem}
Let $({\frak B},{\bf B})$ be an ultrastable quasi--Banach ideal so that $%
{\frak B}\circ {\frak L}_2$ is normed. Let $X_0$ be a Banach space so that $%
X_0$ has the bounded approximation property but $X_0^{\prime \prime }$ has
not. If 
\begin{equation}
{\frak B}(\cdot ,X_0)\subseteq {\frak D}_2(\cdot ,X_0)\text{,}  \tag{$**$}
\end{equation}
then the $({\frak B}\circ {\frak L}_2)^{**}-LRP$ is not satisfied.
\end{theorem}

{\sc Proof:} Put ${\frak A}:=({\frak B}\circ {\frak L}_2)^{*}$. Assume that
the ${\frak A}^{*}-LRP$ is satisfied. Thanks to the previous Lemma, even $%
{\frak A}^{*}\circ {\frak L}_\infty $ has the property (I). Conjugating the
inclusion $(**)$, the total accessibility of ${\frak D}_2$ leads to the
inclusion 
\[
{\frak L}_2(X_0,\cdot )\stackrel{1}{=}{\frak D}_2^\Delta (X_0,\cdot
)\subseteq {\frak B}^\Delta (X_0,\cdot )\stackrel{1}{\subseteq }{\frak B}%
^{*}(X_0,\cdot )\text{,} 
\]
and the quotient formula (\cite{df}, 25.7) therefore implies that $%
Id_{X_0}\in {\frak L}_2^{-1}\circ {\frak B}^{*}(X_0,X_0)\stackrel{1}{=}%
{\frak A}(X_0,X_0)$, and Theorem 3.5 implies the bounded approximation
property of $X_0^{\prime \prime }$ which is a contradiction.$\blacksquare $

Even the case ${\frak B}$ $\subseteq {\frak D}_2^{inj}$ implies the same
situation -- yet requiring a different proof:

\begin{theorem}
Let $({\frak B},{\bf B})$ be an ultrastable quasi--Banach ideal so that $%
{\frak B}$ $\subseteq {\frak D}_2^{inj}$. If ${\frak B}\circ {\frak L}_2$ is
a $1$--Banach ideal, then the $({\frak B}\circ {\frak L}_2)^{**}-LRP$ cannot
be satisfied.
\end{theorem}

{\sc Proof:} As before, put ${\frak A}:=({\frak B}\circ {\frak L}_2)^{*}$.
Assume the validity of the ${\frak A}^{*}-LRP$. Then, ${\frak A}^{*}\circ 
{\frak L}_\infty $ has the property (I), and Theorem 3.4 therefore implies 
\[
({\frak A}^{*\Delta })^{inj}\stackrel{1}{=}({\frak A}^{inj})^{*\Delta }\text{%
.}
\]
On the other hand, since ${\frak B}\subseteq {\frak D}_2^{inj}$, $(*)$
implies 
\[
({\frak B}\circ {\frak L}_2)^{reg}\subseteq {\frak P}_1\text{,}
\]
and it follows 
\[
{\frak L}_\infty \stackrel{1}{=}{\frak P}_1^\Delta \subseteq (({\frak B}%
\circ {\frak L}_2)^{reg})^\Delta \stackrel{1}{=}{\frak A}^{*\Delta }\text{.}
\]
Hence, ${\frak L}\stackrel{1}{=}{\frak L}_\infty ^{inj}\subseteq ({\frak A}%
^{*\Delta })^{inj}\stackrel{1}{=}({\frak A}^{inj})^{*\Delta }\stackrel{1}{%
\subseteq }{\frak I}^\Delta $, and we obtain a contradiction$\blacksquare $.

Permuting the factors ${\frak B}$ and ${\frak L}_2$, we again obtain a
different situation which even shows us a beautyful application of the
principle of local reflexivity for operator ideals to the geometry of Banach
spaces. Let $({\frak B},{\bf B})$ be an ultrastable quasi--Banach ideal so
that ${\frak L}_2\circ {\frak B}$ is normed. Then, we already know that in
this case ${\frak L}_2\circ {\frak B}$ is an injective (hence,
right--accessible) and even maximal Banach ideal, so that the ${\frak L}%
_2\circ {\frak B}-LRP$ automatically is satisfied. Since ${\frak L}_2\circ 
{\frak B}\circ {\frak L}_\infty $ has the property (I), Theorem 3.5 implies
that every Banach space $X$ $\in $ space$({\frak B}^{*}\circ {\frak L}%
_2^{-1})$ with the bounded approximation property even must have a bidual%
{\it \ }$X^{\prime \prime }$ with the bounded approximation property! In
other words:

\begin{theorem}
Let $({\frak B},{\bf B})$ be an ultrastable quasi--Banach ideal so that $%
{\frak L}_2\circ {\frak B}$ is normed. Let $X$ be a Banach space with the
bounded approximation property. If 
\[
{\frak B}(X,\cdot )\subseteq {\frak D}_2(X,\cdot )\text{,}
\]
then even $X^{\prime \prime }$ has the bounded approximation property.
\end{theorem}

To end up this section, we turn to another application of the property (I),
involving Banach spaces of cotype $2$. Using a deep result of Pisier, we
only have to implement some of our own techniques at the right place, to
prove the next result\footnote{%
Note that it is not necessary to assume that $X$ resp. $Y$ has the
Gordon--Lewis property (cf. \cite{djt}, Theorem 17.12).}:

\begin{theorem}
Let $({\frak A},{\bf A})$ be a maximal and left--accessible Banach ideal so
that ${\frak A}^{*}\circ {\frak L}_\infty $ has the property (I). Let $X$
and $Y$ be Banach spaces so that both $X^{\prime }$ and $Y$ have cotype 2.
Then
\end{theorem}

\[
{\frak A}^{inj}(X,Y)\subseteq {\frak L}_2(X,Y)\text{,} 
\]
and 
\[
{\bf L}_2(T)\leq (2{\bf C}_2(X^{\prime })\cdot {\bf C}_2(Y))^{\frac 32}\cdot 
{\bf A}^{inj}(T) 
\]
for all operators $T\in {\frak A}^{inj}(X,Y)$.

{\sc Proof:} Let $X$ and $Y$ be as above and put $C:=(2{\bf C}_2(X^{\prime
})\cdot {\bf C}_2(Y))^{\frac 32}$. Then, \cite[Theorem 4.9.]{pi} tells us,
that any finite rank operator $L\in {\frak F}(Y,X)$ satisfies 
\[
{\bf N}(L)\leq C\cdot {\bf D}_2(L)\text{.}
\]
Hence, 
\[
{\frak N}^\Delta (X,Y)\subseteq {\frak D}_2^\Delta (X,Y)\stackrel{1}{=}%
{\frak L}_2(X,Y)\text{,}
\]
and 
\[
{\bf L}_2(T)\leq C\cdot {\bf N}^\Delta (T)
\]
for all operators $T\in {\frak N}^\Delta (X,Y)$. Since ${\frak N}\stackrel{1%
}{\subseteq }({\frak A}^{inj})^{*}$, we therefore obtain 
\[
({\frak A}^{inj})^{*\Delta }(X,Y)\subseteq {\frak L}_2(X,Y)\text{,}
\]
and 
\[
{\bf L}_2(T)\leq C\cdot ({\bf A}^{inj})^{*\Delta }(T)
\]
for all operators $T\in ({\frak A}^{inj})^{*\Delta }(X,Y)$. Given our
assumptions on ${\frak A}$, Corollary 3.4 reveals that $({\frak A}^{inj})^{*}
$ is totally accessible, and the claim follows.$\blacksquare $

\section{Concluding remarks and open questions}

Summing up our previous investigations, we recognize deep and still
surprising relations between (the validity of) the principle of local
reflexivity for operator ideals, the existence of a norm on product operator
ideals of type ${\frak B}\circ {\frak L}_2$ and the extension of finite rank
operators with respect to a suitable operator ideal norm. The basic objects,
connecting these different aspects, are Jarchow/Otts' product operator
ideals with property (I) and property (S). In the widest sense, a product $%
{\frak A}\circ {\frak B}$ has the property (I), if
\[
({\frak A}\circ {\frak B})\cap {\frak F}=({\frak A}\cap {\frak F})\circ 
{\frak B}
\]
and the property (S), if
\[
({\frak A}\circ {\frak B})\cap {\frak F}={\frak A}\circ ({\frak B}\cap 
{\frak F})\text{,}
\]
so that each finite rank operator in ${\frak A}\circ {\frak B}$ is the
composition of two operators, one of which is of finite rank. Since each
operator ideal which contains ${\frak L}_2$ as a factor, has both, the
property (I) and the property (S), Hilbert space factorization is a
fundamental key. 

However, we do not know whether Corollary 4.2 holds for all maximal Banach
ideals. If this is the case, the ${\frak A}_P^{*}-LRP$ will be false. 

Is the property (I) of ${\frak C}_2^{*}\circ {\frak L}_\infty $ satisfied?
If this is the case, then the injective Banach ideal ${\frak C}_2$ will be
not left--accessible (due to Corollary 2.1 and Corollary 3.4) -- answering
another open question of Defant and Floret (see \cite{df}, 21.2., p. 277).  

We still do not know criteria which are sufficient for the existence of an
ideal--{\it norm} on a given product of quasi--Banach ideals. It seems to be
much easier to give arguments which imply the non--existence of such an
ideal norm (using trace ideals). In particular, we would like to know whether
${\frak T}_2\circ{\frak D}_2$ is a $1$--Banach ideal ($({\frak T}_2, {\bf T}_2)$
denotes the collection of all type 2 operators).


\begin{thebibliography}{99}
\bibitem{a}  {\sf {R. Alencar}}, {\it {Multilinear mappings of nuclear and
intgral type}}, Proc. Amer. Math. Soc. 44 (1985), 33 - 38.

\bibitem{b}  {\sf {A. Braun\ss }}, {\it {Multi-ideals with special
properties }}, Preprint, Potsdam, P\"{a}dagogische Hochschule
''Karl-Liebknecht'' (1987).

\bibitem{cdr}  {\sf {B. Carl, A. Defant, and M. S. Ramanujan}}, {\it {On
tensor stable operator ideals}},\newline
Michigan Math. J. 36 (1989), 63 - 75.

\bibitem{dk}  {\sf D. Dacunha--Castelle and J. L. Krivine}, {\it {\
Applications des ultraproduits \`{a} l'\'{e}tude des espaces et des
alg\`{e}bres de Banach}}, Studia Math. 41 (1972), 315--334.

\bibitem{de}  {\sf {A. Defant}}, {\it {Produkte von Tensornormen}},
Habilitationsschrift, Oldenburg 1986.

\bibitem{df}  {\sf {A. Defant and K. Floret}}, {\it {Tensor norms and
operator ideals}}, North - Holland Amsterdam, London, New York, Tokio 1993.

\bibitem{djt}  {\sf J. Diestel, H. Jarchow, and A. Tonge}, {\it Absolutely
Summing Operators}, Cambridge University Press 1995.

\bibitem{gl}  {\sf {J. E. Gilbert and T. Leih}}, {\it {Factorization, tensor
products and bilinear forms in Banach space theory}}, Notes in Banach
spaces, pp. 182 - 305, Univ. of Texas Press, Austin, 1980.

\bibitem{glr}  {\sf {Y. Gordon, D. R. Lewis, and J. R. Retherford}}, {\it {\
Banach ideals of operators with applications}}, J. Funct. Analysis 14
(1973), 85 - 129.

\bibitem{gr}  {\sf {A. Grothendieck}}, {\it {R\'{e}sum\'{e} de la
th\'{e}orie m\'{e}trique des produits tensoriels topologiques}}, Bol. Soc.
Mat. S\~{a}o Paulo 8 (1956), 1 - 79.

\bibitem{h1}  {\sf {J. Harksen}}, {\it {Tensornormtopologien}},
Dissertation, Kiel 1979.

\bibitem{h}  {\sf {S. Heinrich}}, {\it {Ultraproducts in Banach space theory}
}, J. reine angew. Math. 313 (1980), 72--104.

\bibitem{j}  {\sf {H. Jarchow}}, {\it {Locally convex spaces}}, Teubner 1981.

\bibitem{jo}  {\sf {H. Jarchow and R. Ott}}, {\it {On trace ideals}}, Math.
Nachr. 108 (1982), 23 - 37.

\bibitem{k}  {\sf {K. D. K\"{u}rsten}}, {\it {s--Zahlen und Ultraprodukte
von Operatoren in Banachr\"{a}umen}}, Dissertation, Leipzig, 1976.

\bibitem{l}  {\sf {H. P. Lotz}}, {\it {Grothendieck ideals of operators in
Banach spaces}}, Lecture notes, Univ. Illinois, Urbana, 1973.

\bibitem{lr}  {\sf {J. Lindenstrauss and H. P. Rosenthal}}, {\it {The ${\cal %
L}_p$-spaces}}, Israel J. Math. 7 (1969), 325-349.

\bibitem{oe1}  {\sf {F. Oertel}}, {\it {Konjugierte Operatorenideale und das 
}}${\frak A}${\it {-lokale Reflexivit\"{a}tsprinzip}}, Dissertation,
Kaiserslautern 1990.

\bibitem{oe2}  {\sf {F. Oertel}}, {\it {Operator ideals and the principle of
local reflexivity}}; Acta Universitatis Carolinae - Mathematica et Physica
33, No. 2 (1992), 115 - 120.

\bibitem{oe3}  {\sf {F. Oertel}}, {\it Composition of o{perator ideals and
their regular hulls}}; Acta Universitatis Carolinae - Mathematica et Physica
36, No. 2 (1995), 69 - 72.

\bibitem{oe4}  {\sf {F. Oertel}}, {\it {Local properties of accessible
injective operator ideals}}; Czech. Math. Journal, 48 (123) (1998), 119-133.

\bibitem{p1}  {\sf {A. Pietsch}}, {\it {Operator ideals}}, North - Holland
Amsterdam, London, New York, Tokio 1980.

\bibitem{p2}  {\sf {A. Pietsch}}, {\it {Eigenvalues and s-numbers}},
Cambridge studies in advanced mathematics 13 (1987).

\bibitem{pi}  {\sf {G. Pisier}}, {\it {Factorization of linear operators and
geometry of Banach spaces}}; CBMS Regional Conf. Series 60, Amer. Math. Soc.
1986.
\end{thebibliography}
\end{document}